\documentclass{amsart}

\usepackage{amsthm, amsmath, amssymb, amsbsy, mathdots}
\usepackage{tcolorbox, mathdots}
\usepackage{tikz}
\usepackage{fullpage}
\usepackage{float}
\usepackage[all]{xy}

\newtheorem{theorem}{Theorem}[section]
\newtheorem{lemma}[theorem]{Lemma}

\theoremstyle{definition}
\newtheorem{definition}[theorem]{Definition}

\newtheorem*{theorem*}{Theorem}
\newtheorem*{lemma*}{Lemma}
\newtheorem*{example*}{Example}

\theoremstyle{remark}

\newtheorem{lem}{Lemma}[section]
\newtheorem{cor}[lem]{Corollary}
\newtheorem{rem}[lem]{Remark}

\numberwithin{equation}{section}

\DeclareMathOperator{\GL}{GL}
\DeclareMathOperator{\SL}{SL}
\DeclareMathOperator{\SP}{Sp}
\DeclareMathOperator{\Lie}{Lie}

\renewcommand{\gg}{\mathfrak{g}}
\newcommand{\nn}{\mathfrak{n}}
\renewcommand{\AA}{\mathbb{A}}

\begin{document}

\title{Accessibility of Nilpotent Orbits in Classical Algebras}

\author{Luuk Disselhorst}
\address{}
\curraddr{}
\email{}
\thanks{}


\date{}
\dedicatory{}
\maketitle
\begin{abstract}
Let $G$ be a classical linear algebraic group over an algebraically closed field, and let $\mathfrak{n}$ denote the subset of nilpotent elements in its Lie algebra. In this paper we study a partial order on the $G$-orbits in $\mathfrak{n}$ given by taking limits along cocharacters of $G$. This gives rise to the so-called accessibility order on the nilpotent orbits. Our main results show that for general and special linear algebras, this new order coincides with the usual dominance order on nilpotent orbits, but for symplectic and orthogonal algebras this is not the case.
\end{abstract}

\section{Introduction}
The purpose of this article is to study a new partial order on the nilpotent orbits for classical Lie algebras --- the \emph{accessibility order}. Since the introduction of the classification of nilpotent orbits for algebraic groups, it has been refined often, and this article provides a further step in that refinement. The original classification of Dynkin-Kostant \cite{CollingwoodMcGovern} in characteristic $0$, using $\mathfrak{sl}_2$-triples, is a classic result in the theory, which is still being revisited in the present day, for instance see \cite{ThomasStewart}. In positive characteristic, we refer to the paper of Holt-Spaltenstein (see \cite{Holt-Spaltenstein}) and the sources given there. In this paper, the so-called \emph{dominance order} for nilpotent orbits in classical groups is of particular relevance. The dominance order was first established by Gerstenhaber \cite{Gerstenhaber}.

The idea for using cocharacters (or $1$-parameter subgroups) to study orbits their closures (which can give different results to studying the Zariski-closure in general) finds its origin in the work of Hilbert (where the focus was on the general linear group). The Hilbert-Mumford theorem (see \cite[section 1]{Kempf}) is an application, which says that if an orbit is not Zariski-closed, taking the limit of a suitable cocharacter will yield an element in a different orbit. In the language of this article, this new orbit is then called \emph{1-accessible} from the original orbit. We say that one orbit is \emph{accessible} from another if there exists a chain of cocharacters such that taking limits moves us from one orbit to the other.

Section \ref{Sec:Preliminaries} will cover the preliminaries, regarding the classical groups and algebras, the dominance order, cocharacters and the accessibility order. Section \ref{Sec:ResultsGLn} focuses on the accessibility order of the nilpotent orbits in the general linear algebra, which coincides with the dominance order. Furthermore, in the general linear algebra, 1-accessibility is transitive: that is, if one orbit is accessible from another, then it is 1-accessible. Finally sections \ref{Sec:ResultsSp2n} and \ref{Sec:ResultsOn} cover the results in the symplectic and orthogonal algebras, where the accessibility order and dominance order are not the same, and 1-accessibility is not transitive.
\section{Preliminaries}\label{Sec:Preliminaries}

Our basic references for algebraic groups are Borel \cite{Borel}, Humphreys \cite{Humphreys} and Springer \cite{Springer}.
Let $\kappa$ denote an algebraically closed field, and let $G$ denote a reductive linear algebraic group over $\kappa$; we identify $G$ with its set of $\kappa$-points. 
We say $G$ \emph{acts morphically} on an affine variety $X$ if $G$ acts on $X$ and the action map is a morphism of varieties.
In this case, the point stabilizer $C_G(x)$ for $x\in X$ is a closed subgroup of $G$.

We let $\Lie(G)$ or $\gg$ denote the Lie algebra of $G$, and $\nn$ denote the closed subset of nilpotent elements in $\gg$.
Recall that $G$ acts morphically on $\gg$ and $\nn$ via the adjoint action, and the orbits of $G$ in $\nn$ are called the \emph{nilpotent orbits of $G$}.
In this paper, we are interested in the classical groups of matrices, where this set-up is particularly straightforward to describe; see \cite[p.5]{Jantzen}.

\subsection{General and special linear groups}\label{Prelim:GLnSLn}
If $G = \GL_n(\kappa)$ or $\SL_n(\kappa)$ is the general linear or special linear group in dimension $n$ over $\kappa$, then $\nn$ is the set of nilpotent $n\times n$ matrices, and the action of $G$ is by matrix conjugation.
The nilpotent orbits of $\mathfrak{n}$ are labelled by partitions of $n$, corresponding to the different Jordan Normal forms of nilpotent matrices. When we refer to the orbit of a nilpotent matrix $x$, we denote it with $\mathcal{O}(x)$, or with more emphasis on the action, as $G\cdot x$. Matrices that are in Jordan form will be denoted with a direct sum, e.g. $x=J_{n_1}^{(r_1)}\oplus \cdots \oplus J_{n_p}^{(r_p)}$, where $(r_i)$ is the multiplicity of the standard $n_i \times n_i$ Jordan block $J_{n_i}$. We denote the corresponding partition by $\pi=[n_1^{(r_1)},\ldots,n_p^{(r_p)}]$. A nilpotent matrix induces finite \emph{vector sequences} by acting on a vector. For example, let $x=J_5 \oplus J_3$, then it induces two vector sequences on the standard basis vectors $e_i$, with $1\leq i\leq 8$ (of length 5 and 3, respectively). The sequences are denoted as $x: e_5 \mapsto \ldots \mapsto e_1 \mapsto 0$, and $x: e_8 \mapsto e_7 \mapsto e_6 \mapsto 0$. Comparing induced vector sequences of two matrices, $x$ and $x'$ is a reliable method to check if they are in the same orbit, and to find a conjugating matrix to obtain $x'$ if $x$ is known, see also \cite[Section 2.2]{Disselhorst}. We will be using the vector sequences in Sections \ref{Sec:ResultsGLn}, \ref{Sec:ResultsSp2n} and \ref{Sec:ResultsOn}.

\subsection{Symplectic groups}\label{sec:SymplecticGroups}
We restrict attention to char$(\kappa)\neq 2$.
For concreteness, in what follows we define $G= \SP_{2n}(\kappa)$ to be the group of $2n\times 2n$ matrices satisfying $g^\top \Omega_S g = \Omega_S$, where $\Omega_S$ is the matrix of unit vectors in the following order:  $\Omega_S = (-e_{2n}, \ldots, -e_{n+1}, e_n, \ldots, e_1)$.
Then $\nn$ is the set of nilpotent $n \times n$ matrices satisfying $x^\top\Omega_S+\Omega_S x=0$.
The nilpotent orbits of $G$ are labelled by partitions of $2n$ with the extra condition that each odd part must appear with even multiplicity. A nilpotent element $x\in \mathfrak{g}$ is called \emph{distinguished} if each torus contained in the stabilizer $C_G(x)=\{g\in G \mid g\cdot x = x\}$ is contained in the center of $G$. In the symplectic algebra, $x$ is distinguished if it has a partition $\pi$ where all even parts have multiplicity 1, and no odd parts are present. The process of comparing vector sequences of matrix $x$ and $x'$ will again be used in the symplectic group, but the conjugation obtained from the comparison is not always symplectic. However, elements that are in the same $G$-orbit are also in the same $\text{GL}$-orbit, and it is possible to find a symplectic conjugation from $x$ to $x'$, if any conjugation has been determined. The theorem and process are described in detail in \cite[Sections 1.4-1.5]{Jantzen}.

\subsection{Orthogonal groups}\label{sec:OrthogonalGroups}
Similar to the symplectic groups, we restrict attention to char$(\kappa)\neq 2$. We define $G=O_{n}(\kappa)$ to be the group of $n \times n$ matrices satisfying $g^\top \Omega_O g=\Omega_O$, where $\Omega_O=(e_n,\ldots,e_1)$. Then $\nn$ is the set of nilpotent $n \times n$ matrices satisfying $x^\top\Omega_O+\Omega_O x=0$. In this case, the nilpotent orbits of $G$ are labelled by partitions of $n$ with the extra condition that each even part must appear with even multiplicity. In the orthogonal algebra, $x$ is distinguished if it has a partition $\pi$ where all odd parts have multiplicity 1, and no even parts are present. More information can be found in \cite[1.3-1.6]{Jantzen}.

\subsection{The dominance order}
For any reductive linear algebraic group $G$, there is a natural order on the nilpotent orbits coming from the Zariski topology -- given two orbits $\mathcal{O}_1$ and $\mathcal{O}_1$, we say $\mathcal{O}_2\leq \mathcal{O}_1$ if there is a containment $\overline{\mathcal{O}_2}\subseteq \overline{\mathcal{O}_1}$.
This partial order $\leq$ is called the \emph{dominance order}; it was described in full detail for the classical groups in \cite{Gerstenhaber}.
For $G$ one of the classical groups above, the dominance order on nilpotent orbits matches the dominance order for the partitions labelling the orbits -- that is, if $\mathcal{O}_1$ has partition $\pi_1$, $\mathcal{O}_2$ has partition $\pi_2$, then $\mathcal{O}_2 \leq \mathcal{O}_1$ if and only if $\pi_2 \leq \pi_1$.

\subsection{Cocharacters and accessibility}\label{sec:CocharsAccess}
The principal aim of this paper is to describe \emph{accessibility} between the nilpotent orbits for classical groups, where accessibility is a relation defined by taking limits along cocharacters, and to compare the corresponding partial order to the dominance order. 
We first recall the main details from \cite{BateHerpelMartinRohrle}.
Let $G$ be a reductive linear algebraic group.
A \emph{cocharacter} $\lambda$ of $G$ is a homomorphism of algebraic groups $\lambda:\mathbb{G}_m \to G$, where $\mathbb{G}_m$ denotes the multiplicative group (so $\mathbb{G}_m$ is identified with the nonzero elements of $\kappa$ under multiplication);
we denote the set of cocharacters of $G$ by $Y(G)$.
If $G$ acts morphically on an affine variety $X$, then for every $x\in X$ and $\lambda\in Y(G)$, we have a morphism 
\begin{align*}
\phi_{x,\lambda}:\mathbb{G}_m &\to X\\
a &\mapsto \lambda(a)\cdot x.
\end{align*}
We can now define how $G$ acts on the set of cocharacters. For $g\in G$ and $\lambda \in Y(G)$, we let $g\cdot \lambda$ denote the cocharacter defined by $g\cdot \lambda(a)=g\lambda(a)g^{-1}$, for each $a\in \mathbb{G}_m$
Identifying $\mathbb{G}_m$ as usual as the principal open set in affine $1$-space $\AA^1$ given by throwing away the origin, it makes sense to ask if the morphism $\phi_{x,\lambda}$ extends to a morphism from all of $\AA^1$ to $X$.
When it does, we write $\lim_{a\to0} \lambda(a)\cdot x$ or simply $\lim_\lambda x$ for the image of $0\in \AA^1$ and say \emph{the limit exists}. The Hilbert-Mumford Theorem (see \cite[section 1]{Kempf}) tells us that if an orbit $G\cdot x$ is not closed, there is a cocharacter $\lambda$ such that $\lim_\lambda x \notin G\cdot x$. We can now define accessibility and cocharacter closure.

\begin{definition}
	Suppose $x,y\in X$. $\mathcal{O}(y)$ is called \emph{1-accessible} from $\mathcal{O}(x)$ if there exists a $\lambda \in Y(G)$ such that $\lim_{\lambda}x \in \mathcal{O}(y)$. Furthermore, $\mathcal{O}(y)$ is \emph{$n$-accessible} if there is a chain $x=x_0\rightarrow x_1 \rightarrow \ldots \rightarrow x_n=y$ such that $\mathcal{O}(x_i)$ is 1-accessible from $\mathcal{O}(x_{i-1})$ for $1\leq i \leq n$. If $\mathcal{O}(y)$ is $n-$accessible for some $n$, it is called \emph{accessible}.
\end{definition}

\begin{definition}\label{Def:CocharClosure}
	The cocharacter-closure $$\overline{G\cdot x}^c=\bigcup G\cdot y,$$ where the union is over the orbits $G\cdot y$ that are accessible from $G\cdot x$, is the smallest cocharacter-closed set containing $G\cdot x$.	
\end{definition}

Since accessibility is defined by taking limits, it is clear that if $\mathcal{O}_2$ is accessible from $\mathcal{O}_1$, then $\mathcal{O}_2\leq \mathcal{O}_1$ in the dominance order.
One of the main results of this paper is that the converse is true for the nilpotent orbits of $\GL_n(\kappa)$ and $\SL_n(\kappa)$, see Theorem \ref{thm:OneAccessibleGLn}. This is no longer true for symplectic and orthogonal groups, see Corollary \ref{cor:distnotaccessible} and Remark \ref{rem:distnotaccessible} below.

It is easy to see that accessibility defines a partial order -- it is clearly reflexive (take the trivial cocharacter) and transitive, and antisymmetry follows since if a limit takes us outside a given orbit, then the orbit dimension falls. Note that this argument only works over an algebraically closed field, and whether or not accessibility defines a partial order is an interesting open question over arbitrary fields, see also \cite[section 10]{BateHerpelMartinRohrle}. Finally, note that the relation between Zariski closure and the dominance order is very similar to the relation between cocharacter closure and the accessibility order.

\subsection{Basic facts about cocharacters and limits}
We collect some basic facts. 
The first is the well-known link between cocharacters and parabolic subgroups \cite[p.184]{Humphreys}.

\begin{lem}
Suppose $G$ is a reductive linear algebraic group and $\lambda \in Y(G)$. 
Define the following subsets of $G$:
\begin{align*}
P_\lambda &:= \{g\in G \mid \lim_{a\to 0} \lambda(a)g\lambda(a)^{-1} \textrm{ exists}\},\\
L_\lambda &:= \{g\in G \mid \lim_{a\to 0} \lambda(a)g\lambda(a)^{-1} =g\},\\
U_\lambda &:= \{g\in G \mid \lim_{a\to 0} \lambda(a)g\lambda(a)^{-1}=1\}.
\end{align*}
Then $P_\lambda$ is a parabolic subgroup of $G$, $L_\lambda$ is a Levi subgroup of $P_\lambda$, and $U_\lambda$ is the unipotent radical of $P_\lambda$.
Moreover, for every pair $(P,L)$ consisting of a parabolic subgroup $P$ and a Levi subgroup $L$ of $P$, there exists a cocharacter $\lambda \in Y(G)$ such that 
$=P_\lambda$ and $L = L_\lambda$. 
\end{lem}
The next lemma is used heavily in \cite{BateHerpelMartinRohrle}, for example. For ease of reference, we provide a quick proof.
\begin{lem}\label{lem:Centralizer-OrbitLimit}
Suppose $G$ is a reductive linear algebraic group acting morphically on the affine algebraic variety $X$.
Suppose $x\in X$ and $\lambda\in Y(G)$ are such that $y = \lim_\lambda x$ exists. 
Then $\lambda$ fixes $y$ -- that is $\lambda(a)\cdot y = y$ for all $a\in \mathbb{G}_m$.
In particular, the image of $\lambda$ is a subtorus of the stabilizer $C_G(y)$.
\end{lem}

\begin{proof}
We may embed $X$ $G$-equivariantly inside a vector space $V$ such that the $G$-action is linear (see \cite[Lemma 1.1(a)]{Kempf}). Then V decomposes into a direct sum of weight spaces for $\lambda$. The set of vectors $v\in V$ for which the limit $\lim_\lambda v$ exists is the subspace corresponding to the non-negative weights, and $\lim_\lambda v$ is the component of weight 0. Thus $\lim_\lambda v$ is $\lambda(a)$-fixed for all $a\in \mathbb{G}_m$.
\end{proof}

The next corollary is crucial to our analysis, as it provides the basic obstruction to accessibility for symplectic and orthogonal groups.

\begin{cor}\label{cor:distnotaccessible}
Suppose $G$ is a symplectic or orthogonal group and $y\in \nn$ is a distinguished nilpotent element.
Then it is impossible to find $x\in \nn$ and $\lambda\in Y(G)$ with $\lim_\lambda x = y$. 
\end{cor}

\begin{proof}
If $\lim_\lambda x = y$, then the image of $\lambda$ is a subtorus of $C_G(y)$, by the previous lemma.
But $y$ is distinguished, so $C_G(y)$ contains no non-central subtori. 
\end{proof}

\begin{rem}\label{rem:distnotaccessible}
The corollary above shows immediately that the accessibility order is different to the dominance order for symplectic groups as soon as there are distinguished elements corresponding to partitions other than the top partition $[2n]$.
To see this, note that everything is below this top orbit in the dominance order, but the orbits of other distinguished elements are not accessible from the top orbit.
The first time this becomes apparent for symplectic groups is when $2n = 6$, when we have an orbit of distinguished elements with partition $[4,2]$ which is not accessible from the top orbit with partition $[6]$, even though $[4,2]<[6]$ in the dominance order.

A similar conclusion holds for orthogonal groups.
\end{rem}
Our final preliminary result is also central to the analysis below. 
\subsection{A key technical result}
Suppose that $G$ is a general linear, symplectic or orthogonal group, and $X$ is the Lie algebra. 
By conjugating if necessary, a cocharacter $\lambda$ can be put in any desired standard form -- that is $\lambda$ can be changed to evaluate in the standard diagonal torus, and with powers that are decreasing in size, going down the diagonal.
This means that the parabolic subgroup $P_\lambda$ will be of standard upper block triangular form.
Furthermore, if $x$ is a nilpotent element such that the limit exists, then $x$ is in the Lie algebra of the parabolic subgroup $P_\lambda$. Let $y=\lim_\lambda x$, then since $\lambda$ fixes $y$, we have that $y$ is in the Lie algebra of the Levi subgroup $L_\lambda$ (which is the block diagonal subgroup).

If we take any $g\in L_\lambda$ and replace $x$ with $g\cdot x$, replace $y$ by $g\cdot y$ and $\lambda$ by $g\cdot \lambda$, then if $\lim_{\lambda} x= y$, it follows that $\lim_{g\cdot\lambda} g\cdot x = g\cdot y$. By conjugating with elements from $L_\lambda$, we may further assume that $y$ is in any standard form we choose for nilpotent elements of the Lie algebra of $L_\lambda$. But conjugating $\lambda$ by an element from $L_\lambda$ doesn't change $\lambda$, which means that we can put $y$ in standard form for the Levi subgroup without changing $\lambda$.

Hence we have:
\begin{lemma}
	When considering limits $y=\lim_\lambda x$ in the general linear, symplectic or orthogonal Lie algbras, we may assume that $\lambda$ is diagonal, $y$ is in a standard form, and $x$ is in a corresponding upper block triangular form.
\end{lemma}
We exclusively deal with diagonal cocharacters in classical groups below, so we introduce a shorthand. We write $\lambda = diag(t^{r_1},...,t^{r_s})$ to indicate that $\lambda$ is the cocharacter given by sending $t$ to the given diagonal matrix for each $t\in \mathbb{G}_m$.

\section{Results in the general linear algebra}\label{Sec:ResultsGLn}
The main result of this section is that for orbits in $\mathfrak{gl}_n$, the accessibility order and dominance order are the same, see \ref{thm:OneAccessibleGLn}. Recall that for the general linear group, the nilpotent orbits are labelled by partitions, and the dominance order for the orbits coincides with the usual dominance order for partitions.
\subsection{An example}\label{sec:ExampleGLn}
A key point of all of the work below is that in order to establish accessibility between nilpotent orbits, it is often necessary to conjugate one of the elements away from a standard form. We illustrate this with the example $\mathcal{O}([3,1])\rightarrow \mathcal{O}([2,2])$ in $\mathfrak{gl}_4$, which is not obvious from looking at the standard forms. Let $x'$ be of the form:

$$x'=\left(\begin{array}{cccc}
	0 & 1 & 1 & 0 \\ 
	0 & 0 & 0 & 0 \\ 
	0 & 0 & 0 & 1 \\ 
	0 & 0 & 0 & 0
\end{array}\right).$$
Then it is easy to verify $x'$ is in $\mathcal{O}([3,1])$. The vector sequences $x:e_4 \mapsto e_3 \mapsto e_1 \mapsto 0,$ and $x:e_2-e_1 \mapsto 0$ are of length $3$ and $1$, respectively, hence $x\in \mathcal{O}([3,1])$. Next, if $\lambda=\text{diag}(t,t,1,1)$, then $\lim_\lambda x'\in \mathcal{O}([2,2])$, so $\mathcal{O}([2,2])$ is $1-$accessible from $\mathcal{O}([3,1]).$
\subsection{The general case}
We begin the analysis of accessibility in the general linear algebra by looking at the special case of partitions comprising two parts. The generalization can be made from this first step, and it allows us to set up the notation which will be used for the rest of the paper. Let $$x'=\left(\begin{array}{ccc|ccc}
	&  &  &  &  &  \\
	& J_s &  & 1_{(a)} &  &  \\
	&  &  &  &  &  \\
	\hline
	&  &  &  &  &  \\
	&  &  &  & J_r &  \\
	&  &  &  &  & 
\end{array}\right),$$
where the $J_s$ and $J_r$ are Jordan blocks of size $s$ and $r$, respectably, $(a)$ denotes the location of the nonzero entry in the matrix. In this case $(a)=(k,s+1)$. The matrix $x'$ induces two vector sequences, starting with vector $e_{r+s}$ and $e_s-e_{2s-k}$, respectively, of sizes $r+k$ and $s-k$, respectively. Comparing with the vector sequences of $x=J_{r+1}\oplus J_{s-1}$ allows for determining a base change matrix such that $g\cdot x = x'$, a direct method to show that $x$ and $x'$ are in the same orbit. Finally let $\lambda=(tI_s,I_r)$, then $y=\lim_\lambda x' = J_s\oplus J_r$, hence $y\in \mathcal{O}([r,s])$ and $\mathcal{O}([r,s])$ is $1-$accessible from $\mathcal{O}([r+k,s-k])$.

In the general linear algebra, there are two moves between partitions in the dominance order:

\begin{enumerate}
	\item[I] Multiple parts give pieces to one part, e.g. $[r_1+k_1,r_2+k_2,s-k_1-k_2]\rightarrow [r_1,r_2,s]$.
	\item[II] Multiple parts take pieces from one part, e.g. $[r+k_1+k_2,s_1-k_1,s_2-k_2]\rightarrow [r,s_1,s_2]$.
\end{enumerate}
It can be verified that any move which is not of the above forms, is instead a combination of the two moves. First, it will be shown that move I and II are possible in the accessibility order, i.e. using cocharacters.
\begin{theorem}\label{GLnTwoPiecesToOnePart}
	Let $\pi_1=[\ldots,r_1+k_1,\ldots,r_p+k_p,s-\sum_{i=1}^p k_i,\ldots]$ and $\pi_2=[\ldots,r_1,\ldots,r_p,s,\ldots]$. Then $\mathcal{O}(\pi_2)$ is 1-accessible from $\mathcal{O}(\pi_1)$.
\end{theorem}
\begin{proof}
Let $$x=J_{r_1+k_1}\oplus J_{r_2+k_2} \oplus \cdots \oplus J_{r_p+k_p}\oplus J_{s-\sum_{i=1}^p k_i}.$$ 
And let $$x'=\left(\begin{array}{ccc|cc|c|lcc|l}
	&   &   & 1_{(a_1)} &   &   &   &   &   &   \\
	&J_s&   &   &   & 1_{(b_1)}   & 1_{(c_1)} &   &   &  \\
	&   &   &   &   &   &   &   &   & 1_{(d_1)} \\
	\hline
	&   &   &   &   & -1_{(b_2)}   &   &   &   &\\
	&   &   &J_{r_1} &   &   & -1_{(c_2)} &   &   &\\
	&   &   &   &   &   &   &   &   & -1_{(d_2)}\\
	\hline
	&   &   &   &   &   &   &   &   &\\
	&   &   &   &   & J_{r_2}  & \;\;\vdots &   &   &\\
	&   &   &   &   &   &   &   &   & -1_{(d_3)}\\
	\hline
	&   &   &   &   &   &   &   &   &\\
	&   &   &   &   &   &  \;\;\;\;\;\ddots  &  &   &\vdots\\
	&   &   &   &   &   &   &   &   & \\
	\hline
	&   &   &   &   &   &   &   &   &\\
	&   &   &   &   &   &   &   &   &J_{r_p}\\
	&   &   &   &   &   &   &   &   & 
\end{array}\right),$$
where 

$$\begin{array}{|c|c|}
	\hline 
	(a_1)=(k_1,s+1) & (d_1)=((\sum_{i=1}^pk_i),s+(\sum_{i=1}^{p-1}r_i)+1) \\ 
	\hline 
	(b_1)=(k_1+k_2,s+r_1+1) & (d_2)=(s+(\sum_{i=2}^pk_i),s+(\sum_{i=1}^{p-1}r_i)+1) \\ 
	\hline 
	(b_2)=(s+k_2,s+r_1+1) & (d_3)=(s+r_1+(\sum_{i=3}^pk_i),s+(\sum_{i=1}^{p-1}r_i)+1) \\ 
	\hline 
	(c_1)=((\sum_{i=1}^3k_i),s+(\sum_{i=1}^2r_i)+1) &  \\ 
	\hline 
	(c_2)=(s+k_2+k_3,s+r_1+r_2+1) &  \\ 
	\hline 
\end{array}$$
Then $x'$ induces the following vector sequences:

$$\begin{array}{|c|c|}
\hline
\text{Size} & \text{Starting with} \\
\hline 
r_p+k_p & e_{s+\sum_{i=1}^p r_i} \\
\hline 
r_{p-1}+k_{p-1} & e_{s+\sum_{i=1}^{p-1} r_i} \\
\hline 
\vdots & \vdots \\
\hline
r_1+k_1 & e_{s+r_1} \\
\hline
s-\sum_{i=1}^pk_i & e_s-e_{s+k_1}-e_{s+r_1+k_2}-\ldots-e_{s+(\sum_{i=1}^{p-1})+k_p}\\
\hline
\end{array}$$
Hence $x'\in\mathcal{O}([r_1+k_1,\ldots,r_p+k_p,s-\sum_{i=1}^p k_i])$. Then taking
$$\lambda = \text{diag}(t^pI_s, t^{p-1}I_{r_1}, \ldots, tI_{r_{p-1}}, I_{r+p})$$
gives $y=\lim_\lambda x'\in O([r_1,\ldots,r_p,s])$, so $[r_1,\ldots,r_p,s]$ is 1-accessible from $([r_1+k_1,\ldots,r_p+k_p,s-\sum_{i=1}^p k_i])$.
\end{proof}
\begin{theorem}\label{thm:AccessOnePartToMultiple}
	Let $\pi_1=[\ldots,r+\sum_{i=1}^q k_i,s_1-k_1,\ldots,s_q-k_q,\ldots]$ and let $\pi_2=[\ldots,r,s_1,\ldots,s_q,\ldots]$. Then $\mathcal{O}(\pi_2)$ is 1-accessible from $\mathcal{O}(\pi_1)$.
\end{theorem}
\begin{proof}
Let $x=J_{r+\sum_{i=1}^q k_i} \oplus J_{s_1-k_1} \oplus\ldots \oplus J_{s_q-k_q}$, and let
$$x'=\left(\begin{array}{ccc|ccc|ccc|ccc|cc}
	&   &   &   &   &   &   &   &   &   &   &   &   &     \\
	& J_{s_1} &   & 1_{(a_1)} &   &   &   &   &   &      &   &   &   &  \\
	&   &   &   &   &   &   &   &   &   &   &   &   &     \\
	\hline
	&   &   &   &   &   &   &   &   &   &   &   &   &     \\
	&   &   &   & J_{s_2} &   & 1_{(a_2)}  &   &   &   &   &   &      &  \\
	&   &   &   &   &   &   &   &   &   &   &   &   &     \\
	\hline
	&   &   &   &   &   &   &   &   &   &   &   &   &    \\
	&   &   &   &   &   &   & \ddots &   &   &   &   &   &     \\
	&   &   &   &   &   &   &   &   &   &   &   &   &    \\
	\hline
	&   &   &   &   &   &   &   &   &   &   &   &   &    \\
	&   &   &   &   &   &   &   &   &   & J_{s_q} &   & 1_{(a_q)} &     \\
	&   &   &   &   &   &   &   &   &   &   &   &   &     \\
	\hline
	&   &   &   &   &   &   &   &   &   &   &   &   &     \\
	&   &   &   &   &   &   &   &   &   &   &   &  J_r &    \\
	&   &   &   &   &   &   &   &   &   &   &   &   &     
\end{array}\right),$$
with $(a_1)=(k_1,s_1+1)$, $(a_2)=(s_1+k_2,s_1+s_2+1)$, $((\sum_{i=1}^{q-1}s_i)+k_q,1+\sum_{i=1}^{q}s_i)$.
Then $x'$ induces the following vector sequences:
$$\begin{array}{|c|c|}
	\hline
	\text{Size} & \text{Starting with} \\
	\hline
	r+\sum_{i=1}^q k_i & e_{r+\sum_{i=1}^q s_i} \\
	\hline
	s_q-k_q & e_{\sum_{i=1}^qs_i}-e_{(\sum_{i=1}^q s_i)+s_q-k_q} \\
	\hline
	\vdots & \vdots \\
	\hline
	s_2-k_2 & e_{s_1+s_2}-e_{s_1+2s_2-k_2} \\
	\hline
	s_1-k_1 & e_{s_1+s_2}-e_{2s_1-k_1}\\
	\hline
\end{array}.$$
Hence $x'\in \mathcal{O}([r+\sum_{i=1}^q k_i,s_1-k_1,\ldots,s_q-k_q])$. Then taking
$$\lambda=\text{diag}(t^qI_{s_1},t^{q-1}I_{s_2},\ldots,tI_{s_q},I_{r})$$
gives $y=\lim_\lambda x'\in \mathcal{O}([r,s_1,\ldots,s_q])$, so $[r,s_1,\ldots,s_q]$ is 1-accessible from $[r+\sum_{i=1}^q k_i,s_1-k_1,\ldots,s_q-k_q]$.
\end{proof}
Now the moves of type I and II are verified in the accessibility order, the following theorems show that any move can be described as a combination of the two. Hence, any smaller orbit is 1-accessible from a bigger orbit.
\begin{lemma}\label{Lem:RSTpartition}
	Let $\pi_1=[r+k,s+l-k,t-l]$ and let $\pi_2=[r,s,t]$. Then $\pi_2$ is 1-accessible from $\pi_1$.
\end{lemma}
\begin{proof}First, suppose $k<l$. Then we write $l=k+k'$, and the partition $\pi_1$ can be written as $\pi_1=[r+k,s+k+k'-k ,t-k-k']$, so the move is of type I: $[r+k_1,s+k_2,t-k_1-k_2]\rightarrow [r,s,t]$.
	
Next, suppose $k>l$. Then we write $k=l+k'$ and $\pi_1$ can be written as $\pi_1=[r+l+k,s+l-l-k',t-l]$, so the move is of type II: $[r+k_1+k_2,s-k_1,t-k_2]\rightarrow [r,s,t]$.
	
Finally, if $k=l$, it immediately follows that $s+k-l=s$, so the move is of the form $[r+k,s,t-k]\rightarrow [r,s,t]$, which is clearly a valid move with cocharacters. Hence $\pi_2$ is 1-accessible from $\pi_1$.\end{proof}

\subsection{Conclusion}\label{Sec:GLnConclusion} We can now make several conclusions. First, in $GL_n(\kappa)$, dominance implies $1$-accessibility. Second, $GL_n(\kappa)$-orbits and $SL_n(\kappa)$-orbits are identical. And third, accessibility in $\mathfrak{gl}_n$ implies accessibility in  $\mathfrak{sl}_n$.
\begin{theorem}\label{thm:OneAccessibleGLn}
	Let $\pi_1$ and $\pi_2$ be any two partitions such that $\pi_1$ dominates $\pi_2$. Then $\mathcal{O}(\pi_2)$ is $1$-accessible from $\mathcal{O}(\pi_1)$. Hence, for $GL_n$, accessibility and $1$-accessibility coincide, and the partial order on orbits given by accessibility is the same as the dominance order.
\end{theorem}
\begin{proof}
	$\pi_1$ dominates $\pi_2$, so letting $\pi_1=[a_1,\ldots,a_n]$, $\pi_2=[b_1,\ldots,b_m]$, we have $\sum_{i=1}^n a_i \geq \sum_{i=1}^m b_i$ for $1 \leq i \leq m$ (note also that $n \leq m$, or equivalently, all $a_{n+1},\ldots, a_m$ parts are of size zero).
	
	We can rewrite $\pi_2=[\ldots,r_i,\ldots,s_j,\ldots,t_l,\ldots]$ and $\pi_1=[\ldots,r_i+P_i,\ldots,s_j-Q_j,\ldots,t_l,\ldots]$, so $r_i+P_i$ parts lose pieces, and $s_j-Q_j$ parts gain pieces ($t_l$ parts are unchanged). There are a finite number of $r_i+P_i$ parts, say $p$, and a finite number of $s_j-Q_j$ parts, say $q$.
	
	Then we can denote each $r_i+P_i$ as $r_i+\sum_{j=1}^q p_{i,j}$, where $p_{i,j}$ are the pieces transferred to all $s_j-Q_j$, for $1 \leq j \leq q$, and some $p_{i,j}$ may be zero. Hence we can describe the move $\pi_1 \rightarrow \pi_2$ as $p$ number of moves of type $1$; all $r_i$ parts lose pieces simultaneously. Similarly, all $s_j-Q_j$ parts can be denoted as $s_j-\sum_{i=1}^p q_{i,j}$, hence the move can be described as $q$ moves of type $2$; all $s_j$ parts gain pieces simultaneously. Then we have described all changes to the $r_i$ and the $s_j$ parts in one move, so $\pi_2$ is 1-accessible from $\pi_1$.
\end{proof}
Here follows the accessibility diagram for the nilpotent orbits in $\mathfrak{gl}_6$, that is also the dominance order.
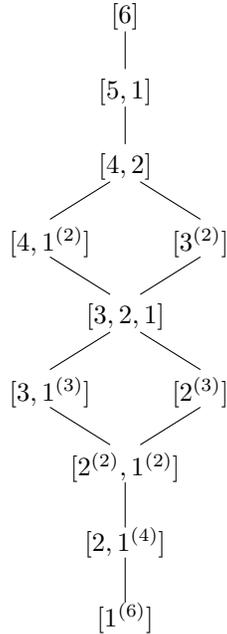
\begin{figure}[ht]
	\centering
	\begin{tikzpicture}
		\node at (2,-1) {$[6]$};
		\draw[-] (2,-1.2) to (2,-1.7);
		\node at (2,-2) {$[5,1]$};
		\draw[-] (2,-2.2) to (2,-2.7);
		\node at (2,-3) {$[4,2]$};
		\draw[-] (1.8,-3.2) to (1,-3.7);
		\draw[-] (2.2,-3.2) to (3,-3.7);
		\node at (1,-4) {$[4,1^{(2)}]$};
		\node at (3,-4) {$[3^{(2)}]$};
		\draw[-] (1,-4.2) to (1.8,-4.7);
		\draw[-] (3,-4.2) to (2.2,-4.7);	
		\node at (2,-5) {$[3,2,1]$};
		\draw[-] (1.8,-5.2) to (1,-5.7);	
		\draw[-] (2.2,-5.2) to (3,-5.7);	
		\node at (1,-6) {$[3,1^{(3)}]$};
		\node at (3,-6) {$[2^{(3)}]$};
		\draw[-] (1,-6.2) to (1.8,-6.7);
		\draw[-] (3,-6.2) to (2.2,-6.7);
		\draw (2,-7.2) to (2,-7.8);
		\node at (2,-7) {$[2^{(2)},1^{(2)}]$};
		\node at (2,-8) {$[2,1^{(4)}]$};
		\draw (2,-8.2) to (2,-8.8);
		\node at (2,-9) {$[1^{(6)}]$};
	\end{tikzpicture}
	\caption{Dominance and $1$-accessibility of partitions in $\mathfrak{gl}_6$.}
	\label{fig:DominanceGL6}
\end{figure}\\
Recall that the general linear group contains as a subgroup the \emph{special linear group} $\text{SL}_n(\kappa)$ consisting of matrices of determinant $1$. The corresponding Lie algebra $\mathfrak{sl}_n(\kappa)$ consists of trace zero matrices in $\mathfrak{gl}_n(\kappa)$. Since the trace of nilpotent matrices is zero, all nilpotent matrices in $\mathfrak{gl}_n$ are also in $\mathfrak{sl}_n$.

The accessibility of nilpotent orbits in the general linear group is helpful for determining the accessibility in the special linear group. In fact, the following two lemmas will show that the accessibility for $\mathfrak{sl}_n$ is the same as that for $\mathfrak{gl}_n$.

\begin{lemma}\label{Lem:GLnSLn-Elements}
	The $\mathfrak{gl}_n(\kappa)$-orbits and $\mathfrak{sl}_n(\kappa)$-orbits of nilpotent matrices are identical.
\end{lemma}
\begin{proof}
	Let $x$ be a nilpotent matrix in standard form, in any orbit, so $x\in \mathfrak{sl}_n \subset \mathfrak{gl}_n$. Since any $g\in \text{GL}_n(\kappa)$ can be written $g = zh$ with $z$ a scalar matrix and $h\in \text{SL}_n(\kappa)$, two nilpotent elements are conjugate by $\text{GL}_n(\kappa)$ if and only if they are conjugate by $\text{SL}_n(\kappa)$.
\end{proof}
Furthermore, the accessibility of orbits in $\mathfrak{sl}_n$ is the same as that for the orbits in $\mathfrak{gl}_n$.
\begin{lemma}\label{Lem:GLnSLn-Access}
	Let $x_1, x_2 \in \mathfrak{gl}_n$ be nilpotent matrices corresponding to partitions $\pi_1$ and $\pi_2$, respectively. If $x_2$ is accessible from $x_1$ in $\mathfrak{gl}_n$, then it is accessible in $\mathfrak{sl}_n$.
\end{lemma}

\begin{proof}
	Let $x_2$ be accessible from $x_1$ in $\mathfrak{gl}_n$, then there is a $g\in \text{GL}_n$ and a cocharacter $\lambda$ of $\text{GL}_n(\kappa)$ such that $\lim_\lambda (g\cdot x_1)=x_2$. From Lemma \ref{Lem:GLnSLn-Elements}, one can find an $h\in \text{SL}_n$ such that $hx_1h^{-1}=gx_1'g^{-1}$. 
	
	The approach to find a suitable $\mu$ is similar. Let $x_2$ have partition $[r_1,\ldots,r_p]$ of $p$ parts, each of size $r_i$. Then $\lambda$ is of the form $\lambda=\text{diag}(t^{p-1}I_{r_{1}},t^{p-2}I_{r_{2}},\ldots,I_{r_p})$, where each $I_{r_{i}}$ is an identity block of size $r_i$. Then 
	
	\begin{align*}
		\det(\lambda)	&=\Pi_{i=1}^{p}\Pi_{j=1}^{r_i} t^{p-i}\\
		&=\Pi_{i=1}^{p}t^{r_i(p-i)}\\
		&=t^{\sum_{i=1}^{p}r_i(p-i)}.
	\end{align*}
	The effect of the cocharacter $\lambda$ when limits are taken does not depend on the values of the powers of $t$, only on the fact that these powers form a decreasing sequence. Hence, if required, we can replace the powers with a decreasing sequence of integers summing to 0 -- that is, we replace each $t^{p-i}$ (including when $i=p$) with powers $t^{a_i}$ with $a_i\in \mathbb{Z}$ chosen so that $\sum_{i=1}^p r_ia_i = 0$. Then the new cocharacter has the same effect in the limit, but now evaluates in $\text{SL}_n(\kappa)$.
	
	Hence, given nilpotent $x_1$ and $x_2$ such that $x_2=\lim_\lambda g\cdot x_1$, one can find $h$ and $\mu$ such that $x_2=\lim_\mu h\cdot x_1$, and $x_2$ is 1-accessible from $x_1$ in $\mathfrak{gl}_n$.
\end{proof}
\subsection{Arbitrary Fields}\label{Sec:GLArbitraryFields}
Our standing assumption is that the field $\kappa$ is algebraically closed. However, this is actually unnecessary for the analysis above of accessibility for $\GL_n$. Let $F$ be an arbitrary field with algebraic closure $\kappa$. We identify $\GL_n(F)$ inside $\GL_n(\kappa)$. Then $\GL_n(F)$ acts by conjugation on nilpotent $n\times n$ matrices with entries in $F$ with exactly the same parametrisation of orbits as in the algebraically closed case. Further, if we restrict attention to $F$-defined cocharacters, we can extend the accessibility order to this setting. The above analysis goes through with minimal changes to show that the accessibility order for $\GL_n(F)$-orbits of nilpotent matrices coincides with the dominance order on the associated partitions. The same is true for the $\SL_n(F)$-orbits.
\section{The symplectic algebra}\label{Sec:ResultsSp2n}
In this section, we analyse the accessibility of nilpotent orbits in the symplectic algebra. First, recall that $\SP_{2n}=\{g\in \GL_{2n}\mid g^\top\Omega_Sg=\Omega_S\}$, where $\Omega_S$ is the matrix of unit vectors in the following order: $\Omega_S=(-e_{2n},\ldots,-e_{n+1},e_n,\ldots,e_1)$. From this point, we will use the notation $f_1=-e_{2n},\ldots,f_n=-e_{n+1}$. The Lie algebra is $\mathfrak{sp}_{2n}=\{x\in\mathfrak{gl}_{2n}\mid x^\top\Omega_S+\Omega_Sx=0\}$. We also recall the possible $\mathfrak{sp}_{2n}$-partitions, and which orbits are distinguished, from Section \ref{sec:SymplecticGroups}.

\begin{theorem}
	Let $x$ be a nilpotent element with partition $[r_1^{(n_1)},\ldots,r_p^{(n_p)}]$. Then $x$ appears in $\mathfrak{sp}_{2n}$ if and only if $n_i$ is even for all odd $r_i$.
\end{theorem}
\begin{lemma}
	Let $x\in \mathfrak{sp}_{2n}$, then $x$ is distinguished if and only if its partition $\pi$ has distinct even parts (and no odd parts).
\end{lemma}
With the following two examples, we show the standard form of the symplectic matrices in this paper. First, let $x$ have partition $\pi=[2n]$, then $x$ is of the form:

$$x=\tilde{J}_{2n}=\left(\begin{array}{ccc|ccc}
	&  &  &  &  &  \\
	& J_n &  &  &  &  \\
	&  &  & 1_{(a)} &  &  \\
	\hline
	&  &  &  &  &  \\
	&  &  &  & \!\!\!\!-J_n &  \\
	&  &  &  &  & 
\end{array}\right),
$$
where $(a)=(n,n+1)$. Next, suppose $x$ has partition $\pi_=[2n_1,2n_2,n_3^{(2)}]$, then $x$ is of the form:
$$x=\left(\begin{array}{c|ccc|ccc|ccc|c}
	J_{n_3} &  &  &  &  &  &  &  &  &  &  \\
	\hline
	&  &  &  &  &  &  &  &  &  &  \\
	&  & J_{n_2} &  &  &  &  &  &  &  &  \\
	&  &  &  &  &  &  & 1_{(a)} &  &  &  \\
	\hline
	&  &  &  &  &  &  &  &  &  &  \\
	&  &  &  &  & \tilde{J}_{2n_1} &  &  &  &  &  \\
	&  &  &  &  &  &  &  &  &  &  \\
	\hline
	&  &  &  &  &  &  &  &  &  &  \\
	&  &  &  &  &  &  &  & \!\!\!\!-J_{n_2} &  &  \\
	&  &  &  &  &  &  &  &  &  &  \\
	\hline
	&  &  &  &  &  &  &  &  &  & -J_{n_3}
\end{array}\right),$$
where $(a)=(n_3+n_2,n_3+n_2+1)$.

Moving between orbits using cocharacters induces corresponding moves between the partitions labelling those orbits. By considering minimal moves in the dominance order, we can work out which are possible to achieve by taking limits along cocharacters.
The five basic moves between adjacent partitions in the dominance order for partitions labelling orbits in the symplectic group are listed below, where we only write down the adjacent parts in the partition which change.
\begin{enumerate}
	\item $\mathcal{O}([2m,2m-2])		\rightarrow \mathcal{O}([2m-1,2m-1])$.
	\item $\mathcal{O}([2n,m,m])		\rightarrow \mathcal{O}([2n-2,m+1,m+1])$.
	\item $\mathcal{O}([n,n,2m])		\rightarrow \mathcal{O}([n-1,n-1,2m+2])$.
	\item $\mathcal{O}([n,n,m,m])		\rightarrow \mathcal{O}([n-1,n-1,m+1,m+1])$.
	\item $\mathcal{O}([2n,2m])			\rightarrow \mathcal{O}([2n-2,2m+2])$.
\end{enumerate}
If the move of type 5 cannot be denoted as a different move (so if $2n \geq 2m+6$), or as a composition of different moves (see  Section \ref{SymplecticNonTransitive}), we call the move a \emph{direct} move of type 5. Otherwise we call it \emph{indirect}. In Section \ref{Subsec:NonMove}, we show that a direct move does not occur. In Section \ref{SymplecticNonTransitive}, we show that while an indirect move is possible, it can only occur in the form of the composition, and not in the form of Move 5 itself (hence 1-accessibility is not transitive).
\subsection{Cocharacters realising moves 1-4}\label{SympPartitionsMoves}
In this section we show that Moves 1-4 above can be achieved using cocharacters. For each of the moves, matrices $x$ and $x'$ will be given, and the starting vectors of the vector sequences that $x'$ induces. Then, the $\lambda$ to obtain $y=\lim_\lambda x'$ will be given. As mentioned in Section \ref{Prelim:GLnSLn}, to show that $x$ and $x'$ are in the same orbit, we will compare the sizes of the vector sequences they induce. Recall from Section \ref{sec:SymplecticGroups} that the conjugating matrix obtained from this comparison is not always symplectic, and the detailed process to find a $g'\in\text{SP}_{2n}$ is described in \cite[Section 1.4-1.5]{Jantzen}.

\emph{Move 1: $\mathcal{O}([2m,2m-2])		\rightarrow \mathcal{O}([2m-1,2m-1])$.} Let:
$$x=\left(\begin{array}{ccc|ccc|cc}
	&  &  &  &  &  &  &    \\ 
	& J_{m} &  &  &  &  &  &    \\ 
	&  &  &  &  &  & 1_{(a)} &    \\ 
	\hline
	&  &  &  &  &  &  &    \\ 
	&  &  &  & \tilde{J}_{2m-2} &  &  &    \\ 
	&  &  &  &  &  &  &    \\ 
	\hline
	&  &  &  &  &  &  &    \\ 
	&  &  &  &  &  & -J_{m} &    \\ 
	&  &  &  &  &  &  &   
\end{array}\right),x'=\left(\begin{array}{ccc|ccc}
	&  &  & 1_{(b)} &  &  \\ 
	& J_{2m-1} &  &  &  &  \\ 
	&  &  &  &  & \!\!\!\!\!\!1_{(c)} \\ 
	\hline
	&  &  &  &  &  \\ 
	&  &  &  & \!\!\!\!\!\!\!\!\!\!-J_{2m-1} &  \\ 
	&  &  &  &  & 
\end{array}\right),$$
with $(a)=(m,3m-1)$, $(b)=(1,2m)$ and $(c)=(2m-1,4m-2)$.
Then $x'$ induces a vector sequence of size $2m$ (starting in $f_1$) and a vector sequence of size $2m-2$ (starting in $f_2-e_n$), hence it is in the orbit $\mathcal{O}([2m,2m-2])$. Taking $\lambda=\text{diag}(tI_{2m-1},t^{-1}I_{2m-1})$ yields $y=\lim_\lambda x'\in \mathcal{O}([2m-1,2m-1])$, hence $\mathcal{O}([2m-1,2m-1])$ is 1-accessible from $\mathcal{O}([2m,2m-2])$.

\emph{Move 2: $\mathcal{O}([2n,m,m])		\rightarrow \mathcal{O}([2n-2,m+1,m+1])$.} Let $x=J_m \oplus\tilde{J}_{2n}\oplus (-J_m)$  and let

$$x'=\left(\begin{array}{ccc|ccc|cc}
	&  &  & 1_{(a)} &  &  &  &   \\ 
	& J_{m+1} &  &  &  &  &  &    \\ 
	&  &  &  &  &  &  &    \\ 
	\hline
	&  &  &  &  &  &  &    \\ 
	&  &  &  & \tilde{J}_{2n-2} &  &  &    \\ 
	&  &  &  &  &  &  & -1_{(b)}   \\ 
	\hline
	&  &  &  &  &  &  &    \\ 
	&  &  &  &  &  & -J_{m+1} &    \\ 
	&  &  &  &  &  &  &   
\end{array}\right),$$
with $(a)=(1,m+2)$ and $b=(m+2n-1,2m+2n)$. Then $x'$ induces three vector sequences, of size $2n$ (starting in $f_1$), $m$ (starting in $f_2$), and $m$ (starting in $e_n-e_{n+m}$). So $x\in \mathcal{O}([2n,m,m])$. Taking $\lambda=\text{diag}(tI_{m+1},I_{2n-2},t^{-1}I_{m+1})$ yields $y=\lim_\lambda x' \in \mathcal{O}([2n-2,m+1,m+1])$, so this orbit is 1-accessible from $\mathcal{O}([2n,m,m])$.

\emph{Move 3:} $\mathcal{O}([n,n,2m])		\rightarrow \mathcal{O}([n-1,n-1,2m+2]).$ Let $x=J_n \oplus \tilde{J}_{2m} \oplus (-J_n),$ and let

$$
x'=\left(\begin{array}{ccc|ccc|cc}
	&  &  &  &  &  &  &  \\ 
	& J_{n-1} &  &  &  &  &  &  \\ 
	&  &  &  &  & \!\!\!\!\!\!\!\!\!1_{(a)} &  &  \\ 
	\hline
	&  &  &  &  &  & 1_{(b)} &  \\ 
	&  &  &  & \;\; \tilde{J}_{2m+2} &  &  &  \\ 
	&  &  &  &  &  &  &  \\ 
	\hline
	&  &  &  &  &  &  &  \\ 
	&  &  &  &  &  & \;\; -J_{n-1} &  \\ 
	&  &  &  &  &  &  & 
\end{array}\right),$$
with $(a)=(n-1,n+2m+1)$ and $(b)=(n,n+2m+2)$. Then $x'$ induces three vector sequences, two of size $n$ (starting in $f_1$, and in $f_n$), and one of size $2m$ (starting in $f_{n+1}+f_{n-2m}$). So $x\in \mathcal{O}([n,n,2m])$. Taking $\lambda=\text{diag}(tI_{n-1},I_{2m+2},t^{-1}I_{n-1})$ yields $y=\lim_\lambda x'\in \mathcal{O}([n-1,n-1,2m+2])$, so this orbit is 1-accessible from $\mathcal{O}([n,n,2m])$.

\emph{Move 4: $\mathcal{O}([n,n,m,m])		\rightarrow \mathcal{O}([n-1,n-1,m+1,m+1])$.} Let $x=J_m\oplus J_n \oplus (-J_n) \oplus (-J_m),$, and let $$x'=\left(\begin{array}{ccc|ccc|ccc|ccc}
	&  &  & 1_{(a)} &  &  &  &  &  &  &  &  \\ 
	& J_{m+1} &  &  &  &  &  &  &  &  &  &  \\ 
	&  &  &  &  &  &  &  &  &  &  &  \\ 
	\hline
	&  &  &  &  &  &  &  &  &  &  &  \\ 
	&  &  &  & J_{n-1} &  &  &  &  &  &  &  \\ 
	&  &  &  &  &  &  &  &  &  &  &  \\ 
	\hline
	&  &  &  &  &  &  &  &  &  &  &  \\ 
	&  &  &  &  &  &  & -J_{n-1} &  &  &  &  \\ 
	&  &  &  &  &  &  &  &  &  &  & 1_{(b)} \\ 
	\hline
	&  &  &  &  &  &  &  &  &  &  &  \\ 
	&  &  &  &  &  &  &  &  &  & -J_{m+1} &  \\ 
	&  &  &  &  &  &  &  &  &  &  & 
\end{array} \right),$$
with $(a)=(1,m+2)$ and $(b)=(m+2n-1,2m+2n).$ Then $x'$ induces four vector sequences, two of size $n$ (starting in $f_1$ and $e_{n+m}$, respectively) and two of size $m$ (starting in $f_2$ and $e_{1+m}-e_{1+2m}$, respectively). So $x\in \mathcal{O}([n,n,m,m])$. Taking $\lambda=\text{diag}(t^2I_{m+1},tI_{n-1},t^{-1}I_{n-1},t^{-2}I_{m+1})$ yields $y=\lim_\lambda x'\in \mathcal{O}([n-1,n-1,m+1,m+1])$, so this orbit is 1-accessible from $\mathcal{O}([n,n,m,m])$.

\subsection{A non-move}\label{Subsec:NonMove}
In this section we analyse Move 5: $\mathcal{O}([2n,2m])\rightarrow \mathcal{O}([2n-2,2m+2])$. We start with examples in section \ref{Subsec:NonMoveExamples}, showing that this move cannot always be achieved with a cocharacter. After that, we set-up some notation to be used in the final section, where we characterise precisely when this move cannot occur.

\subsubsection{Examples}\label{Subsec:NonMoveExamples}First, consider that if $2n-2=2m+2$, then this move is actually a generalized move 1, of the form $\mathcal{O}([2m,2m-2k])\rightarrow \mathcal{O}([2m-k,2m-k])$, with $k=2$. So for move 5, the case where $2n-2\neq 2m+2$ is considered.

Now, consider the move $\mathcal{O}([6])\rightarrow \mathcal{O}([4,2])$ (so $m=0$). From Section \ref{sec:SymplecticGroups}, one can see that the target orbit is distinguished, so it is not accessible from any other orbit, and the move is not valid with cocharacters.
So we can say that in the symplectic algebra, $\overline{G\cdot x}^c \neq \overline{G\cdot x}$.
We finish this example by noting that the presence of extra parts in the partition can complicate matters. For example, even though $\mathcal{O}([4,2])$ is not accessible from $\mathcal{O}([6])$ in $\mathfrak{sp}_6$, we can show that $\mathcal{O}([4,2,2])$ is 1-accessible from $\mathcal{O}([6,2]) \in \mathfrak{sp}_8$. The presence of the extra part of size 2 in the partition makes a material difference to the outcome, which is not obvious immediately. The essential difference in this case is that in this case that the two partitions are no longer adjacent in the dominance order -- the accessibility diagram for this case can be split as follows:
\begin{figure}[ht]
	\centering
	\begin{tikzpicture}
	\node at (1,-1) {$[6,2]$};
	\draw[->] (1.5,-1) to (2.4,-1);
	\draw[->] (1.5,-1.2) to (2.4,-1.7);
	\node at (3,-1) {$[6,1,1]$};
	\draw[->] (3,-1.2) to (3,-1.7);
	\node at (3,-2) {$[4,2,2]$};
	\end{tikzpicture}
	\caption{Decomposition of Move type 5.}
\end{figure}
Consecutively applying move 1 to the $[2]$ part of $[6,2]$ ($[2m,2m-2]\rightarrow [2m-1,2m-1]$ for $m=1$) and move 2 to the partition $[6,1,1]$ ($[2n,m,m] \rightarrow [2n-2,m+1,m+1]$ for $n=3$, $m=1$), shows that $\mathcal{O}([4,2,2])$ is accessible from $\mathcal{O}([6,2])$. Proving 1-accessibility for these orbits is left to the reader. With these examples in mind, a precise analysis of move 5 follows in the next section.

\subsubsection{The Shrinking Operation}\label{Subsec:ShrinkOp}
To prove that move 5 does not occur, a new operation will be introduced. Define the shrinking operation as follows: given a $2n \times 2n$ matrix $A$, matrix $S(A)$ is the $(2n-2) \times (2n-2)$ matrix formed by deleting rows and columns on the outside of matrix $A$. For example, let:

$$A=\left(\begin{array}{cccc}
	a_{(1,1)} & a_{(1,2)} & a_{(1,3)} & a_{(1,4)} \\ 
	a_{(2,1)} & a_{(2,2)} & a_{(2,3)} & a_{(2,4)} \\ 
	a_{(3,1)} & a_{(3,2)} & a_{(3,3)} & a_{(3,4)} \\ 
	a_{(4,1)} & a_{(4,2)} & a_{(4,3)} & a_{(4,4)}
\end{array}\right), \text{ then }
S(A)=\left(\begin{array}{cc}
	a_{(2,2)} & a_{(2,3)} \\ 
	a_{(3,2)} & a_{(3,3)}
\end{array}\right).$$
Denote by $S^d(A)$ the $(2n-2d)\times (2n-2d)$ matrix formed repeating the shrinking operation $d$ times. The shrinking operation is well-behaved with respect to the transpose $S^d(A^T)=(S^d(A))^T$. To show that $S$ is also well-behaved with respect to being in the symplectic Lie algebra, observe that when $\Omega_\text{S}$ is the defining matrix for the form in dimension $2n$, then $S^d(\Omega_\text{S})$ is the defining matrix in dimension $2n-2d$. Then if $A^T\Omega_\text{S}+\Omega_\text{S} A=0,$ it follows that $S^d(A^T)S^d(\Omega_\text{S})+S^d(\Omega_\text{S})S^d(A)=0$. 

Next, suppose $y$ is  in standard form and $\dim(y)=2n$. If $x$ is such that $\lim_{\lambda}x=y$ and $\lambda$ is in standard form, then $\lim_{S^d(\lambda)}S^d(x)=S^d(y)$. There are now two possibilities:
\begin{enumerate}
	\item $S^d(x)$ is conjugate to $S^d(y)$, i.e. no orbit change is made.
	\item $S^d(x)$ is strictly higher than $S^d(y)$ in the dominance order.
\end{enumerate}
With these observations in hand, further reductions can be made.
\subsubsection{The setup}\label{Subsec:setup}
It has already been shown that distinguished partitions are not accessible from any partition higher in the dominance order. The question is: are non-distinguished partitions accessible through a `move' of type 5? This is the case of interest, so the following assumptions can be made:

\begin{enumerate}
	\item[(i)] $y$ is in standard form with a repeated part $[\ldots,d,d,\ldots]$ appearing on the outside of $y$ (so $y$ is not distinguished).
	\item[(ii)] $x$ is another nilpotent element.
	\item[(iii)] $\lambda$ is a diagonal cocharacter in standard form such that $\lim_{\lambda} x=y$.
\end{enumerate}

Under these hypotheses, $x$ can be denoted as $x=y+x_0$, with $\lim_{\lambda} x_0 = 0$. Since $\lambda$ is in standard form, $x_0$ is strictly upper triangular. Furthermore, $\lambda$ must centralize $y$, so $\lambda$ has constant weight $a$ on the first $d$ basis vectors, and constant weight $-a$ on the last $d$ vectors, i.e. for $1\leq i \leq d$, one gets $\lambda(t)e_i=t^ae_i$, $\lambda(t)f_i=t^{-a}f_i$. Let $\lambda_0$ be the cocharacter formed by having weight $a$ on the first $d$ vectors, $-a$ on the last $d$ vectors, and weight zero elsewhere. Then $\lambda_0$ and $\lambda$ are identical on the outside $d$ vectors, and $\lambda_0$ fixes the other vectors.

Then $\lim_{\lambda_0}y=y$ and $x':=\lim_{\lambda_0}x$ exists. Since $\lim_{\lambda}x=y$, $\lim_{\lambda}x'=y$ also, hence $x'$ lies between $x$ and $y$ in the order of dominance. Hence, if $x$ and $y$ are adjacent in the order, there are two possibilities:

\begin{enumerate}
	\item $x'=y$,
	\item $x'$ is conjugate to $x$.
\end{enumerate}

In the case that $x'=y$, the only difference between $x$ and $y$ lies in the outside $d$ rows and columns, i.e. $S^d(x)=S^d(y)$. In the case that $x'$ is conjugate with $x$, one may replace $x$ with $x'$ and assume $x$ has the same repeated $[d,d]$ blocks as $y$ on the outside. The following section will prove that both cases lead to a contradiction, hence move 5 cannot occur in a partition change between adjacent parts.

\subsubsection{The move of type 5}
Here is the key result:
\begin{lemma}\label{lem:Move5Symplectic}
	Let $x\in \mathcal{O}([\ldots,2r,2s,\ldots])$ and let $y \in \mathcal{O}([\ldots,2r-2,2s+2,\ldots])$. Then $\mathcal{O}(y)$ is not accessible from $\mathcal{O}(x)$.
\end{lemma}

\begin{proof}\label{prf:NonMove}Proceed by induction on $n$, where $2n$ is the matrix size. If $n=1$, there are no $r,s$ with $2r\geq 6$. For the inductive step, suppose that for all matrix sizes smaller than $2n$ the claim holds, the following proof is a combination of contradiction and some direct calculation.
	
Suppose there is a $\lambda\in Y(G)$ such that $\lim_{\lambda}x=y$. First note that $y$ is not distinguished, by Corollary \ref{cor:distnotaccessible}. So $y$ must contain a repeated part $[\ldots,d,d,\ldots]$. With the assumption that $y$ and $\lambda$ are in standard form, one may assume that this repeated part appears on the outside of the matrix for $y$, so Section \ref{Subsec:setup} can be used. Since $x$ and $y$ are adjacent, the two cases of that section apply. Starting with the second case, suppose that $x'$ is conjugate to $x$, then $x$ and $y$ share the repeated $[\ldots,d,d\ldots]$ part. With a suitable conjugation, the repeated $d-$parts in $x$ are on the outside as well. Then the shrinking operation applied $d$ times removes that part from both $x$ and $y$. The move between $S^d(x)$ and $S^d(y)$ is of the same form as the move between $x$ and $y$, i.e. it consists of the move $[\ldots, 2r, 2s, \ldots]$ to $[\ldots,2r-2,2s+2,\ldots]$ between adjacent parts, with $2r \geq 2s+6$. Since the matrix sizes have been decreased, this is impossible, by induction. A contradiction follows, hence the second case cannot occur.

Hence case 1 of Section \ref{Subsec:setup} is left. Replacing $x$ with $x'$ and $\lambda$ with $\lambda_0$, one has that $S^d(x)=S^d(y)$ in this case. Note that the first $d$ rows and last $d$ columns of $x$ are related because $x$ is symplectic; if there is an entry $a$ in position $(i,j)$ then $x$ has a related entry in position $(2n+1-j,2n+1-i)$. The main idea is to conjugate $x$ by a suitable symplectic unipotent matrix to kill off most entries in the first $d$ rows which have a further nonzero entry in the column below them, say in row $i$ for $i>d$.

So suppose that $x$ has at least one nonzero entry in column $j$, with $d \leq j \leq 2n-d$, and suppose further that in this column there is another nonzero entry in position $(i,j)$ for some $i>d$. Then, since $x$ looks the same as $y$ apart from the first $d$ rows and last $d$ columns, one can conclude that:

\begin{enumerate}
	\item[(i)] This other entry is the only other nonzero entry in the $j^{\text{th}}$ column, because $S^d(x)=S^d(y)$ is in standard form.
	\item[(ii)] The position of the entry has relation $d < i < j$ because $x$ is strictly upper triangular.
	\item[(iii)] The entry is 1 if $i\leq n$ and $-1$ if $i> n$.
	\item[(iv)] It is the only nonzero entry in the $i^{\text{th}}$ row, except possibly in the last $d$ columns - in particular all entries in row $i$ before the $j^{th}$ column are 0.
\end{enumerate}
To visualise, the matrices $x$ and $y$ are of the following form:
$$\left(\begin{array}{ccc|c|ccc|ccc|ccc|c|ccc}
	&  &  &  & \star & \star & \star & \star & \star & \star & \star & \star & \star &  &  &  &  \\ 
	& J_d &  &  & \star & \star & \star & \star & \star & \star & \star & \star & \star &  &  &  &  \\ 
	&  &  &  & \star & \star & \star & \star & \star & \star & \star & \star & \star &  &  &  &  \\ 
	\hline
	&  &  & \ddots &  &  &  &  &  &  &  &  &  &  &  &  &  \\ 
	\hline
	&  &  &  &  &  &  &  &  &  &  &  &  &  & \star & \star & \star \\ 
	&  &  &  &  & J_{r-1} &  &  &  &  &  &  &  &  & \star & \star & \star \\ 
	&  &  &  &  &  &  &  &  &  & 1_{(a)} &  &  &  & \star & \star & \star \\ 
	\hline
	&  &  &  &  &  &  &  &  &  &  &  &  &  & \star & \star & \star \\ 
	&  &  &  &  &  &  &  & \tilde{J}_{2s+2} &  &  &  &  &  & \star & \star & \star \\ 
	&  &  &  &  &  &  &  &  &  &  &  &  &  & \star & \star & \star \\ 
	\hline
	&  &  &  &  &  &  &  &  &  &  &  &  &  & \star & \star & \star \\ 
	&  &  &  &  &  &  &  &  &  &  & \!\!\!\!\!\!\!\!-J_{r-1} &  &  & \star & \star & \star \\ 
	&  &  &  &  &  &  &  &  &  &  &  &  &  & \star & \star & \star \\ 
	\hline
	&  &  &  &  &  &  &  &  &  &  &  &  & \ddots &  &  &  \\ 
	\hline
	&  &  &  &  &  &  &  &  &  &  &  &  &  &  &  &  \\ 
	&  &  &  &  &  &  &  &  &  &  &  &  &  &  & -J_d &  \\ 
	&  &  &  &  &  &  &  &  &  &  &  &  &  &  &  & 
\end{array}\right).$$
Here the entry $1_{(a)}$ is the entry in the bottom-left location of its respective block, and the stars indicate arbitrary entries at any location in their blocks for $x$, but for $y$ the stars indicate only zero-entries at any location in their blocks. In this situation, let $v$ denote the first $d$ entries in the $j^{\text{th}}$ column, viewed as a column vector of length $d$. Let $v_s$ denote the entry in position $s$ of vector $v$. Let $v_t$ denote the last nonzero entry of $v$, this is position $(t,j)$ in matrix $x$. As $x$ is symplectic, there is a $-v_t$ in position $(2n+1-j,2n+1-t)$. Let $u$ be a unipotent matrix with ones on the diagonal, $-v_t$ in position $(t,i)$, and $v_t$ in position $(2n+1-i,2n+1-t)$. Then the action of $u$ on the basis vectors depends on the position of its non-diagonal entries, hence on the position of $v_t$ and the other nonzero entry in the $j^{\text{th}}$ column. We now split into subcases to cover all possibilities. In each case we define a unipotent element $u$ by describing its action on the basis, with the convention that any basis vector not mentioned is fixed by $u$.

\noindent \emph{Subcase 1, } $j\leq n,$ and $i \leq n$:
\begin{align*}\label{eqn:UnipConj}
	u(e_p)&=
	\begin{cases}
		e_p-v_te_t	&\text{ if } p=i,\\
		e_p			&\text{ otherwise.}
	\end{cases}\\
	u^{-1}(e_p)&=
	\begin{cases}
		e_p+v_te_t	&\text{ if } p=i,\\
		e_p			&\text{ otherwise.}
	\end{cases}
\end{align*}
\emph{Subcase 2, } $j   > n,$ and $i \leq n$:
\begin{align*}
	u(f_p)&=
	\begin{cases}
		f_p+v_tf_i	&\text{ if } p=t,\\
		f_p			&\text{ otherwise.}
	\end{cases}\\
	u^{-1}(f_p)&=
	\begin{cases}
		f_p-v_te_i	&\text{ if } p=t,\\
		f_p			&\text{ otherwise.}
	\end{cases}
\end{align*}
\emph{Subcase 3, } $j   > n,$ and $i    > n$:
\begin{align*}
	u(f_{2n+1-p})&=
	\begin{cases}
		f_{2n+1-p}-v_te_t			&\text{ if } p=i,\\
		f_{2n+1-p}-v_te_{2n+1-i}	&\text{ if } p=2n+1-t,\\
		f_{2n+1-p}					&\text{ otherwise.}
	\end{cases}\\
	u^{-1}(f_{2n+1-p})&=
	\begin{cases}
		f_{2n+1-p}+v_te_t			&\text{ if } p=i,\\
		f_{2n+1-p}+v_te_{2n+1-i}	&\text{ if } p=2n+1-t,\\
		f_{2n+1-p}					&\text{ otherwise.}
	\end{cases}	
\end{align*}
In each of these subcases, when $x$ is conjugated by $u$, the $v_t$ entry is always removed from position $(t,j)$, and a $v_t$ value is always added to the entry in position $(t-1,i)$, so the nonzero value ``moves up and to the left.'' By repeating this conjugation with suitable unipotent matrices, the last nonzero entry in every column $j$ can be removed if this column has a further nonzero entry in position $(i,j)$, hence every nonzero entry in these columns can be removed, apart from the nonzero entry in position $(i,j)$. Then we obtain a new matrix in the same orbit, denoted with $\hat{x}$:
$$\hat{x}=\left(\begin{array}{ccc|c|ccc|ccc|ccc|c|ccc}
	&  &  &  & \star &  &  & \star &  &  &  &  &  &  &  &  &  \\ 
	& J_d &  &  & \star &  &  & \star &  &  &  &  &  &  &  &  &  \\ 
	&  &  &  & \star &  &  & \star &  &  &  &  &  &  &  &  &  \\ 
	\hline
	&  &  & \ddots &  &  &  &  &  &  &  &  &  &  &  &  &  \\ 
	\hline
	&  &  &  &  &  &  &  &  &  &  &  &  &  &  &  &  \\ 
	&  &  &  &  & J_{r-1} &  &  &  &  &  &  &  &  &  &  &  \\ 
	&  &  &  &  &  &  &  &  &  & 1_{(a)} &  &  &  &  &  &  \\ 
	\hline
	&  &  &  &  &  &  &  &  &  &  &  &  &  &  &  &  \\ 
	&  &  &  &  &  &  &  & \tilde{J}_{2s+2} &  &  &  &  &  &  &  &  \\ 
	&  &  &  &  &  &  &  &  &  &  &  &  &  & \star & \star & \star \\ 
	\hline
	&  &  &  &  &  &  &  &  &  &  &  &  &  &  &  &  \\ 
	&  &  &  &  &  &  &  &  &  &  & \!\!\!\!\!\!\!\!-J_{r-1} &  &  &  &  &  \\ 
	&  &  &  &  &  &  &  &  &  &  &  &  &  & \star & \star & \star \\ 
	\hline
	&  &  &  &  &  &  &  &  &  &  &  &  & \ddots &  &  &  \\ 
	\hline
	&  &  &  &  &  &  &  &  &  &  &  &  &  &  &  &  \\ 
	&  &  &  &  &  &  &  &  &  &  &  &  &  &  & -J_d &  \\ 
	&  &  &  &  &  &  &  &  &  &  &  &  &  &  &  & 
\end{array} \right),$$
where $1_{(a)}$ is the entry in the bottom-left of its respective block, and the stars indicate arbitrary entries in the first column or last row of their respective blocks (zero entries are omitted as usual). 

To finish this analysis, we first give an overview and then detail the individual cases which can occur. While $y$ induces vector chains of size $2r-2, 2s+2$ and two chains of size $d$, matrix $\hat{x}$ induces different vector sequences. Consider the effect of $\hat{x}$ on a basis vector $f_i$ for $1 \leq i \leq d$. If $\hat{x} (f_i) \neq 0$, then $f_i$ is in the vector sequence starting with $f_1$, and the value of $f_i$ is added to a vector from another vector sequence in the center of the matrix. In other words, a vector sequence initiated in the $J_d$ Jordan block at the bottom of the matrix will be continued by the $J_{2r-2}$ Jordan blocks or the $J_{2s+2}$ Jordan block in the central portion of the matrix, where in the matrix $y$ the vector chains of parts $[2r-2]$ and $[2s+2]$ are initiated. By the symplectic property, the vector chains of lengths $2r-2$ and $2s+2$ that are terminated in $y$ are picked up by the $J_d$ Jordan block at the top left of the matrix. This shows that matrix $\hat{x}$ has parts $[d-a-b]$ (with multiplicity $2$), $[2r-2+2a]$ and $[2s+2+2b]$ for some $a$ and $b$, contradicting the assumption, hence the move described is not of type 5.

In detail, the following moves occur, in two cases:
\emph{Case 1:} if $s >d$, then $\hat{x}$ induces vector chains of size $s+2a$, $r+2b$ and two chains of size $d-a-b$, with at least one of $a$ and $b$ nonzero. So taking $\lim_\lambda x$ is a move of type 2, in fact two moves of type 2 occur simultaneously:

$$\mathcal{O}([r+2a,s+2b,d-a-b,d-a-b])\rightarrow \mathcal{O}([r,s,d,d]),$$	
is a combination of

\begin{align*}
	\mathcal{O}([r+2a,s+2b,d-a-b,d-a-b]) &\rightarrow \mathcal{O}([r,s+2b,d-b,d-b]),\\
	\mathcal{O}([r,s+2b,d-b,d-b]) &\rightarrow \mathcal{O}([r,s,d,d]).
\end{align*}
\emph{Case 2:} if $d > r$, then there are two subcases, for even $d$ and odd $d$. Without loss of generality, one can assume that the nonzero entry is in the first column after the $[d]$ block, so in the otherwise zero column of the $J_{2r-2}$ block. If $d$ is even, let $d=2d'$, and denote the top left $2d' \times (2d'+1)$ block of matrix $\hat{x}$ as follows:

$$\left(\begin{array}{ccccccccc|c}
	&  &  &  &  &  &  &  &  & x_{1,2d'+1} \\ 
	&  &  &  &  &  &  &  &  & \vdots \\ 
	&  &  &  &  &  &  &  &  & x_{k,2d'+1} \\ 
	&  &  &  &  &  &  &  &  & x_{k+1,2d'+1} \\ 
	&  &  &  & J_d &  &  &  &  & \vdots \\ 
	&  &  &  &  &  &  &  &  & x_{k+i,2d'+1} \\ 
	&  &  &  &  &  &  &  &  & x_{k+i+1,2d'+1} \\ 
	&  &  &  &  &  &  &  &  & \vdots \\ 
	&  &  &  &  &  &  &  &  & x_{k+i+j,2d'+1}\\
	\hline\end{array}\right),$$
We will go through the matrix by row, in reverse order. So we let the nonzero entry be in a different position $\hat{x}_{l,2d+1}$ for each value of $1 \leq l \leq d$ separately, and check what orbit $\hat{x}$ is in. Starting at row $k+i+j$, it will become clear that the move type changes between the nonzero entry being in position $k+i+1$ and $k+i$, and the type changes again between positions $k+1$ and $k$.

When the nonzero entry is in the last row, row $d=2d'$, then $\hat{x}$ has a Jordan block of size $2d'+2d'+2r-2$ (and a Jordan block of size $2s+2$, which is not involved). Then, if one considers the nonzero entry to be in one row higher, this part loses two pieces to two different parts, so $\hat{x}$ has parts of size $2d'+2d'+2r-4$, and there are now two parts of size $1$. For every row higher, the largest part loses two additional pieces to the two smaller parts, until these are of size $2r-2$, hence $\hat{x}$ is in the orbit of $[2d'+2d'-(2r-2),2r-2,2r-2,2s+2]$. So, if the nonzero entry is in one of these last $2r-2$ rows, taking the limit $\lim_\lambda x$ yields the move: 
$$\mathcal{O}([2d'+2d'+(2r-2)-2j,j,j,2s+2]) \rightarrow \mathcal{O}([2d',2d',2r-2,2s+2]),$$
which is a move of type 2. Furthermore, we have determined that $1 \leq j \leq 2r-2$.

Next, if the nonzero entry is another row higher, the first part loses two pieces and the second part gains two pieces, hence taking $\lim_\lambda \hat{x}$ gives the following move: $\mathcal{O}([2d'+2d'-(2r-2)-2,2r-2+2,2r-2,2s+2]) \rightarrow \mathcal{O}([2d',2d',2r-2,2s+2])$. Then, for every the nonzero entry goes up, the first part will be two additional pieces smaller, and the second part will be two additional pieces bigger, hence if the nonzero entry is in row $2d'-(2r-2)-i$, $x$ is in the orbit of $[2d'+2d'-(2r-2)-2i,2r-2+2i,2r-2,2s+2]$, and taking the limit $\lim_\lambda \hat{x}$ yields the move:
$$\mathcal{O}([2d'+2d'-(2r-2)-2i,2r-2+2i,2r-2,2s+2])\rightarrow \mathcal{O}([2d',2d',2r-2,2s+2]),$$
which is a move of type 1. This is possible until $2i=2d'-(2r-2)$, because when this equality holds, $\hat{x}$ is in the orbit of $[2d',2d',2r-2,2s+2]$, which is the same orbit as $y=\lim_\lambda \hat{x}$, and no move occurs. Hence for $1 \leq i \leq d'-(r-1)$, and then for the maximum values of $i$ and $j$, the range of $k$ is as follows: $1\leq k \leq 2d'-i-j = d'-(r-1)$. With this, the following conclusion is reached:

\begin{enumerate}
	\item[(i)] If the nonzero entry is in position $k$, for $1\leq k \leq 2d'-i-j = d'-(r-1)$, then $x$ is in the orbit of $[2d',2d',2r-2,2s+2]$. Taking $\lim_\lambda \hat{x}$ will not change the orbit.
	\item[(ii)] If the nonzero entry is in position $k+i$ for the maximum value of $k$ and $1 \leq i \leq d - (r-1)$, then $\hat{x}$ is in the orbit of $[2d'+2i,2d'-2i,2r-2,2s+2]$. Taking $\lim_\lambda \hat{x}$ is a move of type $1$, where the $[2r-2]$ and $[2s+2]$ pieces are unchanged.
	\item[(iii)] If the nonzero entry is in position $k+i+j$ for maximum values of $k,i$, and $1 \leq j \leq 2r-2$, then $\hat{x}$ is in the orbit of $[2d'+2i+2j,2d'-di-j,2r-2-j]$. Taking $\lim_\lambda \hat{x}$ is a combination of move $1$ and move $2$: $\mathcal{O}([2d'+2i+2j,2d'-di-j,2r-2-j]) \rightarrow \mathcal{O}([2d'+2i,2d'-di,2r-2]) \rightarrow \mathcal{O}([2d,2d,2r-2,2s+2])$.
\end{enumerate}

If $d$ is odd, we denote $d=2d'+1$, and this introduces a row where the nonzero entry places $\hat{x}$ in the orbit of $[2d'+1,2d'-1,r-2,s+2].$ Then the conclusion changes only by the size of the blocks:

\begin{enumerate}
	\item[(i)] If the nonzero entry is in position $k$, for $k \leq d' -(r-1)$, then $\hat{x}$ is in the orbit of $[2d'+1,2d'+1,2r-2,2s+2]$. Taking $\lim_\lambda \hat{x}$ will not change the orbit.
	\item[(ii)] If the nonzero entry is in position $k+i$ for the maximum value of $k$ and $1 \leq i \leq d - (r-1)+1$, then $\hat{x}$ is in the orbit of $[2d'+1+(2(i-1)+1),2d'-(2(i-1)-1),2r-2,2s+2]$. Taking $\lim_\lambda \hat{x}$ is a move of type $1$, where the $[2r-2]$ and $[2s+2]$ blocks are unchanged.
	\item[(iii)] If the nonzero entry is in position $k+i+j$ for maximum values of $k,i$, and $1 \leq j \leq 2r-2$, then $\hat{x}$ is in the orbit of $[2d'+(2(i-1)+1)+2j,2d'-(2(i-1)+1)-j,2r-2-j]$. Taking $\lim_\lambda \hat{x}$ is a combination of move $1$ and move $2$: $\mathcal{O}([2d'+(2(i-1)+1)+2j,2d'-(2(i-1)+1)-j,2r-2-j]) \rightarrow \mathcal{O}([2d'+2i,2d'-di,2r-2]) \rightarrow \mathcal{O}([2d,2d,2r-2,2s+2])$.
\end{enumerate}
Now a conclusion can be drawn regarding the move $\mathcal{O}(\hat{x}) \rightarrow \mathcal{O}(y)$, given there is a a nonzero entry in the $d+1^{\text{th}}$ column: a nonzero entry in the first $d$ rows of the first column above the $J_{r-1}$ Jordan block gives rise to a move involving the parts that are of size $[d]$ and $[2r-2]$ after taking the limit, and we have also established that it does not involve the part of size $[2s+2]$. Similarly, a nonzero entry in the first column above the $\tilde{J}_{2s+2}$ Jordan block gives rise to a move involving the parts that are of size $[d]$ and $[2s+2]$ after taking the limit, but it doesn't involve the part of size $[2r-2]$.

In the case there are nonzero entries in the first columns above both Jordan blocks, taking the limit will give rise to two moves, that separately involve the parts which (after taking the limit) are of size $[2s+2]$ and $[d]$, and the parts which are of size $[2r-2]$ and $[d]$. However, the nonzero entries do not allow for a move between the parts that are of size $[2r-2]$ and $[2s+2]$, after taking the limit. Hence in this special case, the required move does not occur.

So for both case 1 and case 2, the nonzero entries described are only in the first $d$ rows of the matrix so the move only involves the $[2r-2]$ and $[d]$ parts when the nonzero entries are above the $J_{r-1}$ bock. In case the nonzero entries are above the $J_{2s+2}$ block, the move involves the $[2s+2]$ and $[d]$ blocks.

In general, there are $J_{d_i}$ blocks present along the diagonal, for $1 \leq i \leq m$. The nonzero entries in the first column above the $J_{r-1}$ blocks are next to the $J_{d_i}$ blocks, so the move involves the $[2r-2]$ and $[d_i]$ parts (after taking the limit). Regardless of the details of the move, the $[2s+2]$ part is not involved.

Vice versa, the nonzero entries in the first column above the $J_{2s+2}$ Jordan block, in any row of one of the $J_d$ or $J_{d_i}$ Jordan blocks give rise to a move involving the parts $[2s+2],[d],[d_1],\ldots,[d_m]$ (after taking the limit), but the $[2r-2]$ part is not involved.

Hence, regardless of the location of the nonzero entries above the $J_{r-1}$ and $\tilde{J}_{2s+2}$ Jordan block, a move of type $[2r,2s]\rightarrow [2r-2,s2+2]$ does not occur.

This shows that the move $\mathcal{O}([\ldots,r,s,\ldots])\rightarrow \mathcal{O}([\ldots,r-2,s+2,\ldots])$ does not occur; the conclusion is that the orbit $\mathcal{O}([\ldots,r-2,s+2,\ldots])$ is not accessible from the orbit $\mathcal{O}([\ldots,r,s,\ldots])$.
\end{proof}

\subsection{Non-transitivity}\label{SymplecticNonTransitive}
In the previous subsection, we have proven that the move of type 5 does not occur between adjacent partitions $[\ldots,2r,2s,\ldots] \rightarrow [\ldots,2r-2,2s+2,\ldots]$, if $2r \geq 2s+6$. However, if the partitions are not adjacent, the move occurs occurs as a composition two different moves. In the following example, there is a composition of move 2 and move 1. Let $\pi_1=[6,1,1]$, $\pi_2=[4,2,2]$, and $\pi_3=[4,2,1,1]$. Then $\pi_3$ is accessible as follows:

\begin{minipage}{1\textwidth}
	\centering
	\begin{tikzpicture}
	\node at (0.8,-1) {$[6,1,1]$};
	\draw[->] (1.5,-1) to (2.4,-1);
	\node at (3.1,-1) {$[4,2,2]$};
	\draw[->] (3,-1.2) to (3,-1.7);
	\node at (3.1,-2) {$[4,2,1,1]$};
\end{tikzpicture}
\end{minipage}

However, $\mathcal{O}([4,2,1,1])$ is not 1-accessible from $\mathcal{O}([6,1,1])$.
\begin{proof}Let $x\in \mathcal{O}([6,1,1])$ and let $y \in \mathcal{O}([4,2,1,1])$. First, we may conjugate $y$ such that it is in standard form, and we may conjugate $\lambda$ inside $C_G(y)$ such that it is diagonal. Then $y$ is of the form:
$$y = \left(\begin{array}{c|c|cccc|c|c} 
	0 	& 	& 	& 	& 	& 	& 	& \\ \hline
	& 0	& 	&	&	&	& 1	& \\ \hline
	& 	& 0	& 1 	& 	& 	& 	& \\
	& 	& 	& 0 	& 1	& 	& 	& \\
	& 	& 	& 	& 0 	& -1	& 	& \\
	& 	& 	& 	& 	& 0 	& 	& \\ \hline
	& 	& 	& 	& 	& 	& 0 	& \\ \hline
	& 	& 	& 	& 	& 	& 	& 0 
\end{array}\right),$$
and $\lambda=\text{diag}(t^{n},1,1,1,1,1,1,t^{-n})$, with $n\in \mathbb{Z}_{\geq 0}$, fixes $y$. We may assume that $n$ is non-negative (we can swap the order of the 1-blocks). These forms restrict the form of $x$, since we require $\lim_\lambda x = y$, $x$ is of the form:
$$x = \left(\begin{array}{c|c|cccc|c|c} 
	0 	& a	& b	& c	& d	& e	& f	& g		 \\ \hline
	& 0	& 	&	&	&	& 1	& 	f		\\ \hline
	& 	& 0	& 1 	& 	& 	& 	& 	e	\\
	& 	& 	& 0 	& 1	& 	& 	& 	d	\\
	& 	& 	& 	& 0 	& -1	& 	& -c	\\
	& 	& 	& 	& 	& 0 	& 	& 	-b	\\ \hline
	& 	& 	& 	& 	& 	& 0 	& 	-a	\\ \hline
	& 	& 	& 	& 	& 	& 	& 0 
\end{array}\right),$$
with $a,b,c,d,e,f,g \in \mathbb{C}$. Now, we show that no element $x$ of this form can exist in $\mathcal{O}([6,1,1])$. First, $x^5$ has a $b^2$ entry in position $(1,8)$, and zeroes everywhere else, so for a block of size 6 to be present in $x$, we need $b\neq 0$.Then, we row-reduce $x$ to:
$$x \mapsto \left(\begin{array}{c|c|cccc|c|c} 
	0 	& a	& b	& 0	& 0	& 0	& 0	& 0			\\ \hline
		& 0	& 0	&0	&0	&0	& 1	& 	0		\\ \hline
		& 	& 0	& 1 & 0	& 0	& 0	& 	0		\\
		& 	& 	& 0 & 1	& 0	& 0	& 	0		\\
		& 	& 	& 	& 0 & -1& 0	& 0		\\
		& 	& 	& 	& 	& 0 & 0	& 	-b		\\ \hline
		& 	& 	& 	& 	& 	& 0	& 	0		\\ \hline
		& 	& 	& 	& 	& 	& 	& 0 
\end{array}\right).$$
The rows are clearly linearly independent, so the matrix has rank 6. Since an $8 \times 8$ matrix with rank 6 has 2 Jordan blocks, it cannot have partition $[6,1,1]$. Hence $\mathcal{O}([4,2,1,1])$ is not 1-accessible from $\mathcal{O}([6,1,1])$\end{proof}
With this proof, we obtain the following lemma:
\begin{lemma}\label{1-access-nontransitive}
	1-accessibility is not transitive in the symplectic algebras.
\end{lemma}

\subsection{Conclusion for $\mathfrak{sp}_{2n}$}\label{Sec:SympConclusion}
In this section, the five minimal moves in the dominance order of partitions corresponding to symplectic nilpotent orbits have been analysed. These can now describe the order of the orbits in the symplectic algebras. In Section \ref{Sec:GLnConclusion}, the conclusion for accessibility of nilpotent orbits in the general linear group is that the partial order on orbits is the same as the dominance order. This is not the case for the symplectic algebra, as section \ref{Subsec:NonMove} shows that distinguished partitions are not accessible from partitions higher up in the dominance order (e.g. in $\mathfrak{sp}_8$, $\pi_1=[6,2]$ is not accessible from $\pi_2=[8]$). Moreover, if there are two partitions $\pi_1=[\dots,r,s,\ldots]$ and $\pi_2=[\ldots,r-2,s+2,\ldots]$, with $r \geq s+6$, then $\pi_1$ dominates $\pi_2$, but $\pi_2$ is not accessible from $\pi_1$. Since moves 1 -- 4 are valid in the accessibility order (see Section \ref{SympPartitionsMoves}), they give rise to a theorem:

\begin{theorem}\label{Thm:SymplecticAccessibility}
	The partial order on nilpotent symplectic orbits given by the accessibility relation is determined by a combination of the moves of type 1 -- 4.
\end{theorem}
Section \ref{sec:CocharsAccess} tells us that the accessibility order is a partial order (when $\kappa$ is closed), and the accessibility order in the symplectic group can now be related to its dominance order:
\begin{lemma}\label{Lem:SPAccOrderDomOrder}
	Let $x,y$ be two nilpotent symplectic matrices. Then $\overline{\mathcal{O}(y)}^c\leq \overline{\mathcal{O}(y)}^c$ if $\overline{\mathcal{O}(y)}\leq \overline{\mathcal{O}(x)}$, unless the move from $\mathcal{O}(x)$ to $\mathcal{O}(y)$ a direct move of type 5.
\end{lemma}
The difference between the dominance order and the accessibility order of orbits, can be described as follows. Let $O_1$ and $O_2$ be two partitions with $\pi_1=[\ldots,r_1,s_1,\ldots]$ and $\pi_2=[\ldots,r_2,s_2,\ldots]$. Then $O_1 \geq O_2$, while $O_2$ is \emph{not} accessible from $O_1$ if $[\ldots,r_1,s_1,\ldots]=[\ldots,r_2+2k,s_2-2k,\ldots]$ for $k\in \mathbb{Z}_{>0}$, and with $r_1 \geq s_1+6$ and $r_2-2k\neq s_2+2k$.

The following figure shows the accessibility diagram and the dominance order of nilpotent orbits in $\mathfrak{sp}_8$.

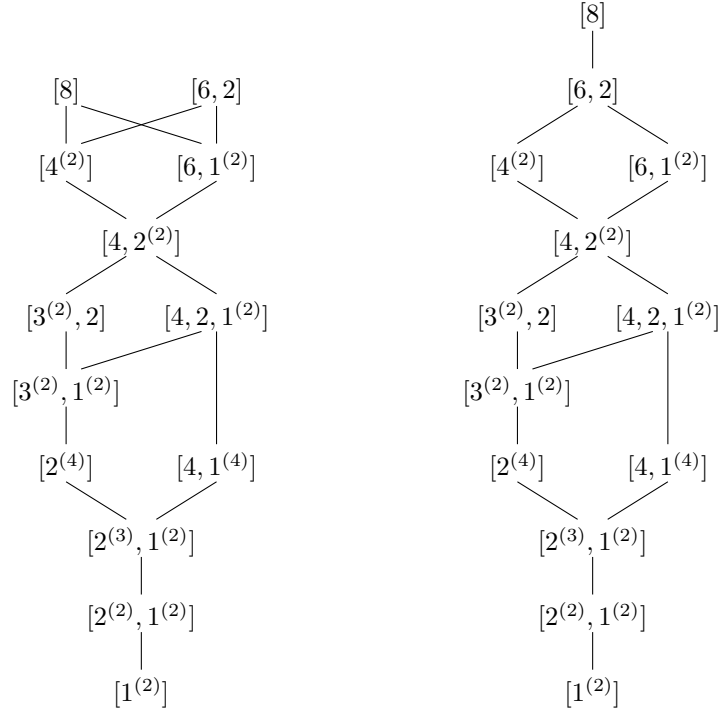
\begin{figure}[ht]
	\centering
	\begin{tikzpicture}
		\node at (1,-1) {$[8]$};
		\node at (3,-1) {$[6,2]$};
		\draw[-] (1,-1.2) to (1,-1.7);
		\draw[-] (3,-1.2) to (3,-1.7);
		\draw[-] (1.2,-1.2) to (2.8,-1.7);
		\draw[-] (2.8,-1.2) to (1.2,-1.7);
		\node at (1,-2) {$[4^{(2)}]$};
		\node at (3,-2) {$[6,1^{(2)}]$};
		\draw[-] (1,-2.2) to (1.8,-2.7);
		\draw[-] (3,-2.2) to (2.2,-2.7);
		\node at (2,-3) {$[4,2^{(2)}]$};
		\draw[-] (1.8,-3.2) to (1,-3.7);
		\draw[-] (2.2,-3.2) to (3,-3.7);
		\node at (1,-4) {$[3^{(2)},2]$};
		\node at (3,-4) {$[4,2,1^{(2)}]$};
		\draw[-] (1,-4.2) to (1,-4.7);
		\draw[-] (2.8,-4.2) to (1.2,-4.7);
		\draw[-] (3,-4.2) to (3,-5.7);	
		\node at (1,-5) {$[3^{(2)},1^{(2)}]$};
		\draw[-] (1,-5.2) to (1,-5.7);
		\node at (1,-6) {$[2^{(4)}]$};
		\node at (3,-6) {$[4,1^{(4)}]$};
		\draw[-] (1,-6.2) to (1.8,-6.7);
		\draw[-] (3,-6.2) to (2.2,-6.7);
		\node at (2,-7) {$[2^{(3)},1^{(2)}]$};	
		\draw[-] (2,-7.2) to (2,-7.7);
		\node at (2,-8) {$[2^{(2)},1^{(2)}]$};
		\draw[-] (2,-8.2) to (2,-8.7);
		\node at (2,-9) {$[1^{(2)}]$};
		
		\node at (8,0) {$[8]$};
		\draw[-] (8,-0.2) to (8,-0.7);
		\node at (8,-1) {$[6,2]$};
		\draw[-] (7.8,-1.2) to (7,-1.7);
		\draw[-] (8.2,-1.2) to (9,-1.7);
		\node at (7,-2) {$[4^{(2)}]$};
		\node at (9,-2) {$[6,1^{(2)}]$};
		\draw[-] (7,-2.2) to (7.8,-2.7);
		\draw[-] (9,-2.2) to (8.2,-2.7);
		\node at (8,-3) {$[4,2^{(2)}]$};
		\draw[-] (7.8,-3.2) to (7,-3.7);
		\draw[-] (8.2,-3.2) to (9,-3.7);
		\node at (7,-4) {$[3^{(2)},2]$};
		\node at (9,-4) {$[4,2,1^{(2)}]$};
		\draw[-] (7,-4.2) to (7,-4.7);
		\draw[-] (8.8,-4.2) to (7.2,-4.7);
		\draw[-] (9,-4.2) to (9,-5.7);	
		\node at (7,-5) {$[3^{(2)},1^{(2)}]$};
		\draw[-] (7,-5.2) to (7,-5.7);
		\node at (7,-6) {$[2^{(4)}]$};
		\node at (9,-6) {$[4,1^{(4)}]$};
		\draw[-] (7,-6.2) to (7.8,-6.7);
		\draw[-] (9,-6.2) to (8.2,-6.7);
		\node at (8,-7) {$[2^{(3)},1^{(2)}]$};	
		\draw[-] (8,-7.2) to (8,-7.7);
		\node at (8,-8) {$[2^{(2)},1^{(2)}]$};
		\draw[-] (8,-8.2) to (8,-8.7);
		\node at (8,-9) {$[1^{(2)}]$};			
	\end{tikzpicture}
	\caption{$1$-accessibility and dominance diagrams of partitions in $\mathfrak{gl}_6$.}
\end{figure}

\section{The orthogonal algebras}\label{Sec:ResultsOn}
\subsection{Partitions and moves}
This section describes some results in the orthogonal algebras. First, recall that $O_N=\{g\in \GL_N \mid g^\top \Omega_O g = \Omega_O\}$, where $\Omega_O$ is the matrix of unit vectors in the following order: $\Omega_O=\{e_N,\ldots,e_1\}$. The Lie algebra is $\mathfrak{o}_N=\{x \in \mathfrak{gl}_N \mid x^\top \Omega_O + \Omega_O x=0\}$. In this section, we will distinguish between even ($N=2n$) and odd ($N=2n+1$) dimension. If $N=2n$, the notation $f_1=e_{2n},\ldots f_n=e_{n+1}$ will be used for the basis vectors. If $N=2n+1$, the notation $f_1=e_{2n+1},\ldots f_n=e_{n+2}$ will be used. We will show that the matrix forms of nilpotent elements are similar to those in the symplectic algebra, and that the orthogonal moves share similarities with the symplectic moves. Note that in this section, we are working over an algebraically closed field $\kappa$, with $\text{char}(\kappa)\neq 2$.

We recall the possible $\mathfrak{o}_N$ partitions, and which orbits are distinguished, from Section \ref{sec:OrthogonalGroups}.
\begin{theorem}
	Let $x$ be a nilpotent element with partition $[r_1^{(n_1)},\ldots,r_p^{(n_p)}]$. Then $x$ appears in $\mathfrak{o}_N$ if and only if $n_i$ is even for all even $r_i$.
\end{theorem}
\begin{lemma}
	Let $x\in \mathfrak{o}_n$, then $x$ is distinguished if and only if its partition $\pi$ has distinct odd parts (and no even parts).
\end{lemma}
As we did with the symplectic algebra, we check two examples to show the standard form of nilpotent orthogonal matrices. First, suppose $x$ has a parition $\pi=[2n+1]$. Then $x$ is of the form:

$$x=\left(\begin{array}{ccccccc}
	0 & 1 &  &  &  &  &  \\
	& \ddots & \ddots &  &  &  &  \\
	&  & \ddots & 1 &  &  &  \\
	&  &  & \ddots & -1 &  &  \\
	&  &  &  & \ddots & \ddots &  \\
	&  &  &  &  & \ddots & -1 \\
	&  &  &  &  &  & 0
\end{array}\right),$$
where the nonzero entries are $-1$ if they are in a row $i$ with $i>\frac{2n+1}{2}$.

Next, let $x$ have partition $\pi=[2n+1,2m+1,2p+1]$. Then $x$ is of the form:

$$x=\left(\begin{array}{ccc|ccc|c|ccc|c|ccc|ccc}
	&  &  &  &  &  &  &  &  &  &  &  &  &  &  &  &  \\
	& J_n &  &  &  &  &  &  &  &  &  &  &  &  &  &  &  \\
	&  &  &  &  &  & \frac{1}{\sqrt{2}} &  &  &  & \frac{1}{\sqrt{2}} &  &  &  &  &  &  \\
	\hline
	&  &  &  &  &  &  &  &  &  &  &  &  &  &  &  &  \\
	&  &  &  & J_m &  &  &  &  &  &  &  &  &  &  &  &  \\
	&  &  &  &  &  & \frac{i}{\sqrt{2}} &  &  &  & -\frac{1}{\sqrt{2}} &  &  &  &  &  &  \\
	\hline
	&  &  &  &  &  &  &  &  &  &  & \frac{i}{\sqrt{2}} &  &  & -\frac{1}{\sqrt{2}} &  &  \\
	\hline
	&  &  &  &  &  &  &  &  &  &  &  &  &  &  &  &  \\
	&  &  &  &  &  &  &  & \tilde{J}_{2p+1} &  &  &  &  &  &  &  &  \\
	&  &  &  &  &  &  &  &  &  &  &  &  &  &  &  &  \\
	\hline
	&  &  &  &  &  &  &  &  &  &  & -\frac{i}{\sqrt{2}} &  &  & -\frac{1}{\sqrt{2}} &  &  \\
	\hline
	&  &  &  &  &  &  &  &  &  &  &  &  &  &  &  &  \\
	&  &  &  &  &  &  &  &  &  &  &  & -J_m &  &  &  &  \\
	&  &  &  &  &  &  &  &  &  &  &  &  &  &  &  &  \\
	\hline
	&  &  &  &  &  &  &  &  &  &  &  &  &  &  &  &  \\
	&  &  &  &  &  &  &  &  &  &  &  &  &  &  & -J_n &  \\
	&  &  &  &  &  &  &  &  &  &  &  &  &  &  &  & 
\end{array}\right),$$
where the $\pm \frac{1}{\sqrt{2}}$ and $\pm\frac{i}{\sqrt{2}}$ are in the last rows to the right of the $J_n$ and $J_m$ blocks, and the first columns above the $-J_n$ and $-J_m$ blocks.

Like in section \ref{SympPartitionsMoves}, here follow the minimal moves through the dominance order of the partitions corresponding to orthogonal nilpotent orbits:
\begin{enumerate}
	\item $\mathcal{O}([2m+1,2m-1])		\rightarrow \mathcal{O}([2m,2m])$.
	\item $\mathcal{O}([2n+1,m,m])		\rightarrow \mathcal{O}([2n-1,m+1,m+1])$.
	\item $\mathcal{O}([n,n,2m-1])		\rightarrow \mathcal{O}([n-1,n-1,2m+1])$.
	\item $\mathcal{O}([n,n,m,m])		\rightarrow \mathcal{O}([n-1,n-1,m+1,m+1])$.
	\item $\mathcal{O}([2n+1,2m+1])		\rightarrow \mathcal{O}([2n-1,2m+3])$.
\end{enumerate}

Comparing with the moves in Section \ref{SympPartitionsMoves}, it shows that the moves in the orthogonal algebra are similar to those in the symplectic algebra; only the sizes of the partitions that occur once are different. In moves 2, 3, and 4 the matrix forms are very similar as well, therefore these moves are omitted. In the next section, Move 1 is analysed, because the standard form of a matrix with partition $[2m+1,2m-1]$ is quite different from the standard form of a matrix with partition $[2m,2m-2]$, and it is worth viewing the differences in more detail. Move 5 is impossible to realise with cocharacters, like with the symplectic algebra, this is a non-move. Like in move 1, the difference in the proof for the non-move is due to the different standard form. For the proof, see \cite[section 4.2]{Disselhorst}.

\subsection{The cocharacter realizing move 1}\label{SubSec:OrthMove1}
In first move for the orthogonal algebras is: $\mathcal{O}([2m+1,2m-1]) \rightarrow \mathcal{O}([2m,2m])$. We start with example $\mathcal{O}([7,5])\rightarrow \mathcal{O}([6,6])$. The standard form of an element in $\mathcal{O}([7,5])$ is:

$$x=\left(\begin{array}{cccccccccccc}
	0 & 1 & 0 & 0 & 0 & 0 & 0 & 0 & 0 & 0 & 0 & 0 \\
	0 & 0 & 1 & 0 & 0 & 0 & 0 & 0 & 0 & 0 & 0 & 0 \\
	0 & 0 & 0 & 0 & 0 & \frac{1}{\sqrt{2}} & \frac{1}{\sqrt{2}} & 0 & 0 & 0 & 0 & 0 \\
	0 & 0 & 0 & 0 & 1 & 0 & 0 & 0 & 0 & 0 & 0 & 0 \\
	0 & 0 & 0 & 0 & 0 & \frac{i}{\sqrt{2}} & -\frac{i}{\sqrt{2}} & 0 & 0 & 0 & 0 & 0 \\
	0 & 0 & 0 & 0 & 0 & 0 & 0 & \frac{i}{\sqrt{2}} & 0 & -\frac{1}{\sqrt{2}} & 0 & 0 \\
	0 & 0 & 0 & 0 & 0 & 0 & 0 & -\frac{i}{\sqrt{2}} & 0 & -\frac{1}{\sqrt{2}} & 0 & 0 \\
	0 & 0 & 0 & 0 & 0 & 0 & 0 & 0 & -1 & 0 & 0 & 0 \\
	0 & 0 & 0 & 0 & 0 & 0 & 0 & 0 & 0 & 0 & 0 & 0 \\
	0 & 0 & 0 & 0 & 0 & 0 & 0 & 0 & 0 & 0 & -1 & 0 \\
	0 & 0 & 0 & 0 & 0 & 0 & 0 & 0 & 0 & 0 & 0 & -1 \\
	0 & 0 & 0 & 0 & 0 & 0 & 0 & 0 & 0 & 0 & 0 & 0 \\
\end{array}\right).$$
The matrix $x$ induces the following vector sequences.
\begin{align*}
	f_1 \rightarrow -f_2 \rightarrow f_3 \rightarrow &\frac{1}{\sqrt{2}}(-f_6-e_6) \rightarrow -e_3 \rightarrow -e_2 \rightarrow -e_1 \rightarrow 0,\\
	f_4 \rightarrow -f_5 \rightarrow &\frac{i}{\sqrt{2}}(f_6-e_6) \rightarrow e_5 \rightarrow e_4 \rightarrow 0.
\end{align*}
Next, let $x'$ be of the form:

$$x'=\left(
\begin{array}{cccccccccccc}
	0 & 1 & 0 & 0 & 0 & 0 & 1 & 0 & 0 & 0 & 0 & 0 \\
	0 & 0 & 1 & 0 & 0 & 0 & 0 & 0 & 0 & 0 & 0 & 0 \\
	0 & 0 & 0 & 1 & 0 & 0 & 0 & 0 & 0 & 0 & 0 & 0 \\
	0 & 0 & 0 & 0 & 1 & 0 & 0 & 0 & 0 & 0 & 0 & 0 \\
	0 & 0 & 0 & 0 & 0 & 1 & 0 & 0 & 0 & 0 & 0 & 0 \\
	0 & 0 & 0 & 0 & 0 & 0 & 0 & 0 & 0 & 0 & 0 & -1 \\
	0 & 0 & 0 & 0 & 0 & 0 & 0 & -1 & 0 & 0 & 0 & 0 \\
	0 & 0 & 0 & 0 & 0 & 0 & 0 & 0 & -1 & 0 & 0 & 0 \\
	0 & 0 & 0 & 0 & 0 & 0 & 0 & 0 & 0 & -1 & 0 & 0 \\
	0 & 0 & 0 & 0 & 0 & 0 & 0 & 0 & 0 & 0 & -1 & 0 \\
	0 & 0 & 0 & 0 & 0 & 0 & 0 & 0 & 0 & 0 & 0 & -1 \\
	0 & 0 & 0 & 0 & 0 & 0 & 0 & 0 & 0 & 0 & 0 & 0 \\
\end{array}
\right),$$
then it induces the vector sequences:
\begin{align*}
	f_1 \rightarrow -f_2-e_6 \rightarrow f_3 -e_5 \rightarrow -f_4-e_4 &\rightarrow f_5-e_3 \rightarrow -f_6-e_2 \rightarrow -2e_1 \rightarrow 0.\\
	f_2 -e_6 \rightarrow -f_3-e_5 \rightarrow f_4-e_4 &\rightarrow f_5-e_3 \rightarrow f_6-e_2 \rightarrow 0.
\end{align*}
Hence $x\in \mathcal{O}([7,5])$. Let $\lambda=\text{diag}(t,\ldots,t,t^{-1},\ldots,t^{-1})$, then $y=\lim_\lambda x'\in\mathcal{O}([6,6])$, so $\mathcal{O}([6,6])$ is 1-accessible from $\mathcal{O}([7,5])$. We leave the proof for the general case to the reader.
\subsection{Conclusion for $\mathfrak{o}_n$}\label{OrthConclusion}
As we have observed the many similarities between accessibility of nilpotent orbits in the symplectic algebra, and the accessibility of nilpotent orbits in the orthogonal algebra, we can draw the same conclusion. Theorem \ref{Thm:SymplecticAccessibility} in Section \ref{Sec:SympConclusion} applies to the orthogonal algebra as well:

\begin{theorem}\label{Thm:OrthogonalAccessibility}
	The partial order on nilpotent orthogonal orbits given by the accessibility relation is determined by a combination of the moves of type 1 -- 4.
\end{theorem}
\begin{proof}
	In the previous section, we have seen that move 1 is a valid move in the orthogonal algebra, and section \ref{SympPartitionsMoves} suffices to show that moves 2, 3 and 4 are valid in the orthogonal algebra, due to the similarities with the respective symplectic moves. Finally, the proof that move 5 is impossible to realise in $\mathfrak{o}_n$, is also very similar to the proof of its symplectic counterpart.
\end{proof}
Since the accessibility order is a partial order (when $\kappa$ is closed), it can be related to the dominance order.
\begin{lemma}\label{Lem:OrthAccOrderDomOrder}
	Let $x,y$ be two nilpotent orthogonal matrices. Then $\overline{\mathcal{O}(y)}^c\leq \overline{\mathcal{O}(y)}^c$ if $\overline{\mathcal{O}(y)}\leq \overline{\mathcal{O}(x)}$, unless the move from $\mathcal{O}(x)$ to $\mathcal{O}(y)$ is of type 5.
\end{lemma}

\section{Extensions and open questions}
In this paper we have analysed the accessibility order of nilpotent orbits in the classical algebras. In the general linear and special linear algebras, the accessibility order coincides with the dominance order of partitions, and 1-accessibility is transitive. In the symplectic and orthogonal group, the accessibility order differs from the dominance order, and 1-accessibility is not transitive. The work is restricted to algebraically closed fields with characteristic not 2. When these restrictions are lifted, there are open questions about accessibility in the symplectic and orthogonal groups.

\subsection{Non-algebraically closed fields}
The work is restricted to algebraically closed fields $\kappa$, with $\text{char}(\kappa)\neq 2$. There are open questions about the extent to which the assumption of an algebraically closed field can be relaxed. This is important, since one of the motivations for the introduction of cocharacter closure is to provide a formalism to work over an arbitrary field. For example, studying Zariski-closed sets may fail to pick up interesting behaviour, as described in \cite[Example 3.4]{BateHerpelMartinRohrle}.

We begin with a reminder: in Section \ref{Sec:GLArbitraryFields}, we have shown that the accessibility order of nilpotent orbits in the general and special linear algebras are the same when the restriction of closed fields is lifted. 

The situation is less transparent for symplectic and orthogonal groups, as one of the key results - the theorem (see \cite[Section 1.4]{Jantzen}) that two elements $x$ and $x'$ in one $G$-orbit are also in a $\text{GL}$-orbit, may fail over non-closed fields. The process of finding a symplectic (or orthogonal) conjugation from $x$ to $x'$ may not work, even though a $\text{GL}$-conjugation is known. In \cite[Section 5]{Disselhorst}, the symplectic conjugation determine for the move $\mathcal{O}([6,4]) \rightarrow \mathcal{O}([5,5])$, which involves $\sqrt{2}$ and $\sqrt{-1}$, that are not always defined over non-closed $\kappa$. Furthermore, the standard forms of matrices in the orthogonal algebras involve entries with $\sqrt{2}$ and $\sqrt{-1}$ values, which gives an extra challenge if the restriction of close fields is lifted.

This means that we cannot analyse nilpotent orbits over an arbitrary field without additional tools at our disposal. However, there are some options if we restrict our attention to certain classes of fields.

As an example, further research could be fruitful in the direction of finite fields. Let $G$ be the symplectic or orthogonal group defined over $\kappa=\overline{\mathbb{F}}_q$, the closure of the finite field with $q$ elements. Then $\mathbb{F}_q$ can be realised as the fixed points of the Frobenius endomorphism $x\mapsto x^q$, which extends to a homomorphism $\sigma$ which raises each matrix entry to the $q^{\text{th}}$ power:
\begin{align*}\sigma: G &\rightarrow G,\\
	\text{given by }a_{ij}	&\mapsto a_{ij}^q.
\end{align*}

It is clear that the subgroup $G_\sigma$ of $\sigma$-fixed points in $G$ is the finite symplectic or orthogonal group consisting of matrices satisfying the defining conditions for $G$, but with all entries in $\mathbb{F}_q$. We can now extend the Frobenius endomorphism to a map on the Lie algebra $\mathfrak{g}$ of $G$ and consider how the $G-$orbits of $\mathfrak{g}$ relate to the $G_\sigma-$orbits of $\mathfrak{g}_\sigma$.

In particular, if $x\in \mathfrak{g}$ is such that $\sigma(x)=x$, how does $(G\cdot x)_\sigma:=G\cdot x \cap \mathfrak{g}_{\sigma}$ decompose into $G_\sigma-$orbits? The theorem of Springer-Steinberg (see \cite[p.172-173]{Springer-Steinberg}) gives the answer:

\begin{theorem*}
	Let $C=C_G(x)$ be the centralizer of $x$ in $G$. Then the $G_\sigma-$orbits in $(G\cdot x)_\sigma$ are parametrised by the elements of the cohomology group $H^1(\sigma,C)$. In the special case that $C$ is connected, there is a single orbit.
\end{theorem*}
Here $H^1(\sigma,C)$ denotes $C$ modulo the equivalence relation $a \sim b$ if and only if there exists $c\in C$ with $a=cb\sigma(c)^{-1}$. These are the orbits under `twisted' conjugation. With this result and the knowledge of the centralizers of nilpotent elements from \cite[Section 3]{Jantzen}, in principle we can work out the nilpotent orbits for these finite groups and then begin to analyse accessibility relations between them. In practice, this is likely to be a complex process, but it is possible that we could find a general theorem analogous to theorem \ref{Thm:SymplecticAccessibility}, which deals with accessibility of nilpotent symplectic orbits over a closed field (or the likely similar theorem for nilpotent orthogonal orbits).
\bibliographystyle{amsplain}

\end{document}